\documentclass{amsart}
\usepackage{amsmath, amssymb, amsthm, gensymb}
\usepackage{hyperref}
\usepackage{pgfplots, tikz, tikz-cd}
\usepackage{xfrac, relsize}
\usepackage{footnote}
\usepackage{subcaption}
\usepackage{enumitem}
\usepackage[dvipsnames]{xcolor}
\usepackage[backend=biber, citetracker=true, sorting=nty, sortcites=true]{biblatex}

\calclayout

\usepackage[lf]{Baskervaldx}              
\usepackage[bigdelims, vvarbb]{newtxmath} 
\usepackage[cal=boondoxo]{mathalfa}       

\setlist[enumerate,1]{label={\normalfont (\arabic*)}}

\hypersetup{
  colorlinks,
  citecolor = blue,
  urlcolor  = blue,
  linkcolor = blue,
}

\pgfplotsset{compat=1.16}
\usetikzlibrary{calc, arrows}

\newcommand{\Diff}{\operatorname{Diff}^+}
\DeclareMathOperator{\Mod}{Mod}
\DeclareMathOperator{\PMod}{PMod}

\DeclareMathOperator{\SIP}{SIP}

\DeclareMathOperator{\GL}{GL}

\DeclareMathOperator{\Sp}{Sp}

\DeclareMathOperator{\range}{range}
\DeclareMathOperator{\rank}{rank}

\renewcommand{\le}{\leqslant}

\renewcommand{\ge}{\geqslant}

\renewcommand{\subset}{\subseteq}

\renewcommand{\chi}{\ensuremath \raisebox{\depth}{$\mathchar"11F$}}

\newcommand{\CC}{\mathbb{C}} 
\newcommand{\DD}{\mathbb{D}} 
\newcommand{\RR}{\mathbb{R}} 
\newcommand{\ZZ}{\mathbb{Z}} 

\newtheorem{theorem}{Theorem}[section]
\newtheorem{main-theorem}{Theorem}
\newtheorem{main-corollary}[main-theorem]{Corollary}
\newtheorem{proposition}[theorem]{Proposition}

\newtheorem{lemma}[theorem]{Lemma}
\newtheorem{corollary}[theorem]{Corollary}

\theoremstyle{remark}
\newtheorem*{remark}{Remark}

\numberwithin{equation}{section} 

\newcommand{\torushole}[2]{
  \begin{scope}[shift={(#1,#2)}]
    \draw[thick, line cap=round]
      (210:.4)+(0,.25)
      arc[x radius=.4, y radius=.3, start angle=210, end angle=330];
    \draw[thick] (135:.3)+(0,-.27) arc (135:45:.3);
  \end{scope}
}

\title{Lower Bounds for Faithful Linear Representations of Subgroups of the Mapping Class Group}
\author{Thiago Brevidelli}
\address{Universit\'e de Toulouse, CNRS, Institut de Math\'ematiques de Toulouse, France}
\email{thiago.brevidelli\_garcia@math.univ-toulouse.fr }
\urladdr{https://www.math.univ-toulouse.fr/~tbrevide/}
\subjclass[2020]{Primary 57K20; Secondary 20F38 15A30}

\begin{document}

\maketitle

\begin{abstract}
  Recently, Korkmaz established the lower bound of \(3g - 2\) for the dimension
  of a faithful representation of the mapping class group of an orientable
  surface of genus \(g \ge 3\). We raise this bound to \(4g - 3\) in the
  setting of surfaces of genus \(g \ge 7\). A new ingredient is a finer study
  of the commutation relations in $\PMod(\Sigma)$. We use the relations
  arising from a certain pants decomposition of $\Sigma_g$ to show that any
  representation of dimension $\le 4g - 4$ is forced to kill a natural subgroup
  of the Torelli group.

  We also establish lower bounds for the dimension of faithful representations
  of related groups: the Johnson group of a closed surface, arbitrarily low
  terms of the Johnson filtration of a compact surface with one boundary
  component and pure braid groups. These lower bounds grow linearly on the
  genus of the surfaces and the number of strands of the braids. Finally, we
  also provide some evidence that greater lower bounds for the low-genus cases
  should lead to improved lower bounds for \(g \gg 0\).
\end{abstract}
\keywords{mapping class group, representation theory}

\tableofcontents

\section{Introduction}

Denote by \(\Sigma = \Sigma_{g, r}^b\) the compact connected orientable surface
of genus \(g\) with \(b\) boundary components and \(r\) marked points \(P =
\{x_1, \ldots, x_r\}\) in its interior. Its \emph{mapping class group} is the
group \(\Mod(\Sigma) = \pi_0 \Diff(\Sigma, \partial \Sigma)_P\) of
orientation-preserving self-diffeomorphisms of \(\Sigma\) fixing the boundary
point-wise and permuting the marked points, up to isotopies. The \emph{pure
mapping class group} \(\PMod(\Sigma) \le \Mod(\Sigma)\) is the subgroup of
mapping classes fixing \(P\) point-wise. Let \(\Sigma_g = \Sigma_{g, 0}^0\).

The theory of mapping class groups plays a central role in low-dimensional
topology, as closed \(3\)-manifolds may be encoded by elements of
\(\Mod(\Sigma_g)\) via mapping tori or Heegaard splittings. The group
\(\PMod(\Sigma)\) is also of importance in algebraic geometry, as it may be
seen as the (orbifold) fundamental group of the moduli space of complex curves.

Yet, basic questions about its linear representations \(\PMod(\Sigma) \to
\GL_d(\CC)\) remain unanswered. Recent years have seen intense activity
around the study of low-dimensional representations of \(\PMod(\Sigma)\), as
well as the question of linearity of such groups -- which is currently open for
\(g \ge 3\).

Improving results of Funar \cite{funar-11} and Franks--Handel
\cite{franks-handel-13}, Korkmaz \cite{korkmaz-23} showed any \(d\)-dimensional
linear representations of \(\PMod(\Sigma)\) with \(d < 2g\) is trivial for \(g
\ge 3\). Korkmaz then went on to show that, for \(g \ge 3\), any nontrivial
\(\PMod(\Sigma) \to \GL_{2g}(\CC)\) is conjugate to the so-called
\emph{symplectic representation} \(\Psi : \PMod(\Sigma) \to \Sp_{2g}(\ZZ)\).
When \(\Sigma\) has no marked points, its kernel is a normal subgroup of
primary importance, known as \emph{the Torelli subgroup \(\mathcal{I}(\Sigma)
\le \PMod(\Sigma)\)}.

Denoting by \(d(G)\) the smallest dimension of a faithful linear representation
of a group \(G\) and setting \(d(\Sigma) = d(\PMod(\Sigma))\), Korkmaz also
established \(d(\Sigma) \ge 3g - 2\) for \(g \ge 3\). He showed that, when \(g
\ge 3\) and \(m \le g - 3\), any \(\PMod(\Sigma) \to \GL_{2g + m}(\CC)\) is
forced to kill the \(m^{\text{th}}\) derived subgroup \(K_{\Sigma'}^{(m)}\) of
a certain subgroup \(K_{\Sigma'} \le \mathcal{I}(\Sigma)\), where \(\Sigma'
\subset \Sigma\) is a genus \(3\) subsurface.

Kasahara \cite{kasahara-23} classified all \((2g+1)\)-dimensional
representations for \(g \ge 7\). Recently, Kaufmann--Salter--Zhang--Zhong
\cite{kaufmann-25} further improved Korkmaz' results by classifying all
\(\PMod(\Sigma) \to \GL_d(\CC)\) with \(d \le 3g - 3\) in the \(g \ge 4\) and
\(b + r \le 1\) setting, showing that any such representation is conjugate to
the direct sum of a \((2g+1)\)-dimensional representation with copies of the
trivial representation \(\PMod(\Sigma) \to \GL_1(\CC)\).

Their result shed light into Korkmaz' lower bound of \(3g - 2\) in the setting
of closed unmarked surfaces, showing any \(\Mod(\Sigma_g) \to \GL_d(\CC)\) with
\(d \le 3g - 3\) is forced to kill the entire Torelli subgroup.

In this article we raise Korkmaz' lower bound to \(4g - 3\) in the setting of
surfaces of genus \(g \ge 7\). We show that, when \(d\) is small enough, any
\(\PMod(\Sigma) \to \GL_d(\CC)\) is forced to kill the subgroup
\(\SIP_0(\Sigma) \le \mathcal{I}(\Sigma)\) generated by the commutators \([T_a,
T_b]\) of Dehn twists \(T_a, T_b \in \PMod(\Sigma)\) about pairs of curves \(a,
b \subset \Sigma\) intersecting at two points, with algebraic intersection
pairing \(\langle a, b \rangle = 0\) and \(\Sigma \setminus (a \cup b)\)
connected.

\begin{main-theorem}[Theorem~\ref{thm:bound-on-faithful-mcg-real}]\label{thm:bound-on-faithful-mcg}
  Let \(\Sigma\) be a surface of genus \(g \ge 7\) and \(\rho : \PMod(\Sigma)
  \to \GL_d(\CC)\). If \(d \le 4g-4\) then \(\SIP_0(\Sigma) \le \ker \rho\). In
  particular, \(d(\Sigma) \ge 4g - 3\).
\end{main-theorem}

\begin{remark}
  The assumption of \(g \ge 7\) is only used to handle the \(d = 4g - 4\) case.
  As a consequence, the same statement holds if we take \(g \ge 4\) and \(d <
  4g - 4\) -- see Theorem~\ref{thm:bound-on-faithful-mcg-real}. The author
  believes Theorem~\ref{thm:bound-on-faithful-mcg} should hold as stated for
  \(g \ge 4\).
\end{remark}

Unlike the subgroups \(K_{\Sigma'}^{(m)}\) from Korkmaz' proof, the subgroup
\(\SIP_0(\Sigma) \le \ker \rho\) from our proof remains the same regardless of
\(d\). It is a natural subgroup of the group \(\SIP(\Sigma) \le
\mathcal{I}(\Sigma)\) generated by the so-called \emph{simple intersection
maps}: the commutators of twists about curves intersecting at two points and
whose algebraic intersection number vanishes.

These maps were introduced by Putman in \cite{putman-07} as part of a
generating set for the Torelli group of an unmarked genus \(0\) surface. Putman
would then go on to use such maps in his infinite presentation of
\(\mathcal{I}(\Sigma_g)\) \cite{putman-09}. The groups \(\SIP(\Sigma_g)\) and
\(\SIP(\Sigma_g^1)\) were also investigated in their own right by Childers
\cite{childers-12}, who proposed a systematic study of their properties.

We also establish lower bounds for the dimensions of faithful representations
of related groups. Our proofs are elementary in nature, relying mostly on well
known facts about surface mapping class groups.

A new ingredient is a finer study of the commutation relations in
\(\Mod(\Sigma)\). We make use of such relations and certain families of curves
to produce quotients of \(F_2 \times \cdots \times F_2\), the direct product of
\(n\) copies of a rank-\(2\) free group, inside different subgroups \(G \le
\Mod(\Sigma)\).

Such quotients are then used to bound the dimensions of faithful
representations of certain \(G \le \Mod(\Sigma)\) by the smallest dimension of
a faithful representation of \(F_2 \times \cdots \times F_2\). The latter was
recently computed by Kionke--Schesler \cite{kionke-schesler-23}, who showed
that the dimension of a faithful representation of \(F_2 \times \cdots \times
F_2\) is \(\ge 2n\).

Taking the following subgroups for convenience and denoting
\(
  d(G) = \min \{ d | \rho : G \hookrightarrow \GL_d(\CC) \text{ is faithful} \}
\)
as above, we arrive at the theorems
bellow.
\begin{enumerate}
  \item The \emph{Johnson subgroup} \(\mathcal{K}(\Sigma_g) \le
    \mathcal{I}(\Sigma_g)\).
  \item The terms
    \(
      \mathcal{I}(\Sigma_g^1) = \mathcal{I}^0(\Sigma_g^1)
      \triangleright \mathcal{I}^1(\Sigma_g^1)
      \triangleright \cdots
      \triangleright \mathcal{I}^k(\Sigma_g^1)
      \triangleright \cdots
    \)
    of the \emph{Johnson filtration}.
  \item The \emph{pure} braid group \(PB_n\) on \(n\) strands.
\end{enumerate}

\begin{main-theorem}[Corollary~\ref{thm:small-dim-rep-kills-separating-twists}]\label{thm:bound-on-faithful-johnson}
  Let \(g \ge 2\). Then \(d(\mathcal{K}(\Sigma_g)) \ge 2g - 2\).
\end{main-theorem}

\begin{main-theorem}[Corollary~\ref{thm:bound-on-faithful-johnson-filtration-real}]\label{thm:bound-on-faithful-johnson-filtration}
  Let \(g \ge 2\) and \(k \ge 1\). Then \(d(\mathcal{I}^k(\Sigma_g^1)) \ge 2g -
  2\).
\end{main-theorem}

\begin{main-theorem}[Corollary~\ref{thm:bound-on-faithful-braid-real}]\label{thm:bound-on-faithful-braid}
  If \(n\) is odd then \(d(PB_n) \ge n - 1\). If \(n\) is even then \(d(PB_n)
  \ge n - 2\).
\end{main-theorem}

\begin{remark}
  The lower bound of \(n - 1\) from Theorem~\ref{thm:bound-on-faithful-braid}
  is well known for the full braid group \(B_n\). Dyer--Formanek--Grossman
  showed in \cite[Proposition~2]{dyer-82} that if \(\rho : B_n \hookrightarrow
  \GL_d(\CC)\) is faithful then one of its irreducible subquotients must also
  be faithful. The irreducible representations of \(B_n\) of dimension at most
  \(n - 1\) were classified by Formanek \cite{formanek-96}. In particular, no
  irreducible representation of dimension \(< n - 1\) is faithful for \(n \ge
  3\).
\end{remark}

Using a similar strategy, we provide some evidence that greater lower bounds
for the low-genus cases should lead to improvements of the lower bounds in
Theorem~\ref{thm:bound-on-faithful-mcg}.

\begin{main-theorem}[Theorem~\ref{thm:bigger-bounds-real}]\label{thm:bigger-bounds}
  Let \(n \ge 1\) and \(g \ge 2n\). Then \(d(\Sigma_g^1) \ge n \cdot \min \{
  d(E): E \text{ is a cyclic extension of } \Mod(\Sigma_{\lfloor \sfrac{g}{n}
  \rfloor, 1}) \}\).
\end{main-theorem}

\begin{remark}
  To the best of the authors knowledge, the value of \(d(\Sigma_{g, 1})\) is
  unknown even for \(g = 2\). The group \(\Mod(\Sigma_2^1)\) contains a natural
  copy of \(B_5\), so that \(d(\Sigma_2^1) \ge d(B_5)\). The smallest faithful
  representation of \(B_5\) known in the literature seems to be the so-called
  \emph{Lawrence representation} \(B_5 \hookrightarrow \GL_{10}(\ZZ[q^{\pm 1},
  t^{\pm 1}])\) \cite{bigelow-00, krammer-02}. If one assumed \(d(E) \ge 10\)
  for all other cyclic extensions \(E\) of \(\Mod(\Sigma_{2,1})\) then
  Theorem~\ref{thm:bigger-bounds} would say \(d(\Sigma_g^1) \ge 5g\), thus
  improving our lower bounds of \(4g - 3\).
\end{remark}

In fact, the above method is quite general, and one could imagine applying it
to many other groups of interest in low-dimensional topology. Taking \(G =
\PMod(\Sigma)\) for a surface \(\Sigma\) of even genus \(g \ge 4\) even, we
recover Korkmaz' lower bound of \(3g - 2\).

To go beyond this bound we instead use a different strategy. We consider a
family of simple closed curves \(a_1, \ldots, a_{3g-3}, b_1, \ldots, b_{3g-3}
\subset \Sigma\) where the curves \(a_i\) come from a certain pants
decomposition of \(\Sigma_g\), while the curves \(b_j\) are ``complementary''
to the curves \(a_i\) -- see \textsection\ref{sec:main-result} for a
definition.

Given \(\rho : \PMod(\Sigma) \to \GL_d(\CC)\) with \(d\) small enough, we show
that, unless \(\SIP_0(\Sigma) \le \ker \rho\), the matrices \(M_i =
\rho(T_{a_i}) - 1\) and \(N_j = \rho(T_{b_j}) - 1\) satisfy the ``annihilation
relations''
\begin{align}\label{eq:nm-rels-intro}
  N_j M_i & = 0 \iff i \ne j    &
  M_j N_i & = 0 \iff i \ne j    &
  M_j M_i & = 0 \; \forall i, j,
\end{align}
where \(T_a\) denotes the Dehn twist about \(a \subset \Sigma\). This is
accomplished using the disjointness relations in \(\PMod(\Sigma)\), as well as
a careful study of the eigenvalues and eigenspaces of \(\rho(T_a)\) following
the work of Korkmaz, Kasahara and Kaufmann--Salter--Zhang--Zhong.

We establish lower bounds for \(d\) such that we can find \(M_1, \ldots, M_n,
N_1, \ldots, N_n \in M_d(\CC)\) satisfying (\ref{eq:nm-rels-intro}). Together
with the previous assertion about \(M_i = \rho(T_{a_i}) - 1\) and \(N_j =
\rho(T_{b_j}) - 1\), such lower bounds show \(\rho\) is indeed forced to kill
\(\SIP_0(\Sigma)\), thus concluding the proof.

\subsection{Outline of the Paper}

In \textsection\ref{sec:background} we review the theory of mapping class
groups needed for the rest of the paper. This includes some cohomological
calculations, used to handle the \(d = 4g - 4\) case of
Theorem~\ref{thm:bound-on-faithful-mcg}. The informed reader is invited to
skip this section entirely if so inclined.

In \textsection\ref{sec:matrix-graphs} we study commutation relations in
\(\PMod(\Sigma)\). We use these relations and certain families of curves to
produce quotients of \(F_2 \times \cdots \times F_2\) inside of different \(G
\le \PMod(\Sigma)\). We then use these subgroups to establish the lower bounds
from Theorem~\ref{thm:bound-on-faithful-johnson},
Theorem~\ref{thm:bound-on-faithful-johnson-filtration} and
Theorem~\ref{thm:bound-on-faithful-braid}. We then adapt this strategy to
establish Theorem~\ref{thm:bigger-bounds}.

Still in \textsection\ref{sec:matrix-graphs}, we establish a lower bound for
\(d\) such that we can find \(M_1, \ldots, M_n, N_1, \ldots, N_n \in M_d(\CC)\)
satisfying (\ref{eq:nm-rels-intro}) (Proposition~\ref{thm:dim-inf-bound}). In
\textsection\ref{sec:eigenvals} we study the eigenspaces of \(L_a =
\rho(T_a)\), where \(a \subset \Sigma\) is nonseparating and \(\rho :
\PMod(\Sigma) \to \GL_d(\CC)\) is low-dimensional. We establish a lower bound
for the dimension of the \(1\)-eigenspace of \(L_a\)
(Proposition~\ref{thm:eigenval-multiplicity-quota}).

Finally, in \textsection\ref{sec:main-result} we conclude our proof of
Theorem~\ref{thm:bound-on-faithful-mcg} by applying
Proposition~\ref{thm:dim-inf-bound} and
Proposition~\ref{thm:eigenval-multiplicity-quota}.

\subsection*{Acknowledgments}

I would like to express my deep gratitude to my advisors J. Marché and
M. Wolff for taking me as their student and suggesting me this problem, as
well as for the incredibly helpful discussions we shared in the past year.

I should also thank F. Costantino for suggesting me to consider representations
of the Torelli group. I am thankful to B. Petri for communicating me the
Kaufmann--Salter--Zhang--Zhong preprint, and to G. Massuyeau for communicating
me Childers' results on the groups \(\SIP(\Sigma_g)\) and \(\SIP(\Sigma_g^1)\).
Finally, I am grateful to J. Martel and the anonymous referees for their
helpful suggestions on earlier drafts.

This work received support from the University Research School EUR-MINT (State
support managed by the National Research Agency for Future Investments program
bearing the reference ANR-18-EURE-0023). This work was also supported by the
MathPhDInFrance International Doctoral Training Program.

\section{Background Results on Mapping Class Groups}
\label{sec:background}

Let \(\Sigma = \Sigma_{g, r}^b\) be the compact surface of genus \(g\) with
\(b\) boundary components and \(r\) marked points in its interior. We denote by
\(P\) its set of marked points. We may freely omit \(b\) and \(r\) from the
notation when they are zero. All surfaces considered in the present paper have
the form \(\Sigma_{g, r}^b\).

The mapping class group \(\Mod(\Sigma) = \pi_0 \Diff(\Sigma, \partial
\Sigma)_P\) is the group of orientation-preserving self-diffeomorphisms of
\(\Sigma\) up to isotopy, where both our diffeomorphisms and isotopies are
assumed to fix the boundary point-wise and permute the marked points. The group
\(\Mod(\Sigma)\) acts on the set \(P\), and the pure mapping class group
\(\PMod(\Sigma)\) is the subgroup of mapping classes acting trivially.

The \emph{braid group} on \(n\) strands is \(B_n = \Mod(\DD_n)\), the mapping
class group of a disk \(\DD_n = \Sigma_{0, n}^1\) with \(n\) marked points.
This is isomorphic to the the fundamental group of the configuration space of
\(n\) unordered points in a disk. The \emph{pure braid group} on \(n\) strands
is \(PB_n = \PMod(\DD_n)\).

In this section we collect the results from the theory of mapping class groups
needed in the rest of the paper. We refer the reader to \cite{farb-margalit}
for further information on mapping class groups.

Given a closed subsurface \(\Sigma' \subset \Sigma\) with marked points \(P'
\subset P\), there is an induced group homomorphism \(\iota : \PMod(\Sigma')
\to \PMod(\Sigma)\). Such a homomorphism needs not be injective, but we
nevertheless refer to the post-composition of \(\rho : \PMod(\Sigma) \to
\GL_d(\CC)\) by \(\iota\) as \emph{the restriction of \(\rho\) to
\(\PMod(\Sigma')\)}.

\subsection{Curves \& Dehn Twists}
\label{sec:dehn-twists}

Given an unoriented simple closed curve \(\alpha \subset \Sigma\) avoiding the
marked points of \(\Sigma\), we denote by \(a\) its free homotopy class and
write ``\(a \subset \Sigma\)''. Here our homotopies are assumed to
avoid the marked points. All curves considered in this paper are simple closed
curves, unless explicitly stated otherwise.

Recall that the \emph{geometric intersection number} between \(a, b \subset
\Sigma\) is the infimum
\[
  |a \pitchfork b|
  = \min \{ \# (\alpha \pitchfork \beta) : \alpha \in a, \beta \in b\}
\]
of the number of times two transverse representatives of \(a\) and \(b\) cross
each other. On the other hand, the \emph{algebraic intersection number}
\(\langle a, b \rangle\) between \(a\) and \(b\) is the sum of the indices of
the intersection points \(x \in \alpha \pitchfork \beta\) for any \(\alpha \in
a\) and \(\beta \in b\).

We denote by \(T_a\) the \emph{right Dehn twist about \(a\)}. This is the class
of a diffeomorphism of \(\Sigma\) supported in an annular neighborhood of
\(\alpha \in a\) which ``winds a full turn around \(\alpha\).''

It is also useful consider the twists about curves parallel to the boundary
components of \(\Sigma\). For example, by collapsing the boundary into a marked
point we obtain a surjective group homomorphism \(\Mod(\Sigma_g^1) \to
\Mod(\Sigma_{g, 1})\). Its kernel is the subgroup generated by the Dehn twist
\(T_d\) about the boundary \(d = \partial \Sigma_g^1\).

Improving results of Hatcher--Thurston \cite{hatcher-thurston-80} and Harer
\cite{harer-83}, Wajnryb \cite{wajnryb-83, birman-wajnryb-94} produced a
remarkable finite presentation of \(\PMod(\Sigma)\), whose generators are given
by Dehn twists and whose relations can all be explained in terms of the topology
of \(\Sigma\). In this paper we only need a small fragment of this result.

\begin{theorem}[Dehn--Lickorish, \cite{dehn-38, lickorish-62, lickorish-64}]\label{thm:mcg-generators}
  The group \(\PMod(\Sigma)\) is generated by finitely many Dehn twists about
  nonseparating simple closed curves.
\end{theorem}

We also summarize some of the most useful relations in \(\PMod(\Sigma)\).

\begin{enumerate}
  \item \textbf{The conjugation relation.} Given \(a \subset \Sigma\) and
    \(f \in \PMod(\Sigma)\), \(T_{f(a)} = f T_a f^{-1}\).
  \item \textbf{The disjointness relation.} Given \(a, b \subset \Sigma\),
    \(T_a\) commutes with \(T_b\) if and only if \(|a \pitchfork b| = 0\),
    i.e. if and only if we can find disjoint representatives for \(a\) and
    \(b\).
  \item \textbf{The braid relation.} Given \(a, b \subset \Sigma\) with
    \(|a \pitchfork b| = 1\), \(T_a T_b T_a = T_b T_a T_b\).
\end{enumerate}

Here a crucial observation is due: given nonseparating \(a, b \subset \Sigma\),
\(T_a\) and \(T_b\) are conjugate in \(\PMod(\Sigma)\). Indeed, we can always
find \(f \in \PMod(\Sigma)\) such that \(f(a) = b\), so that \(f T_a f^{-1} =
T_b\) by the conjugation relation. Together with
Theorem~\ref{thm:mcg-generators}, this implies the Abelianization
\(\PMod(\Sigma)^{\operatorname{ab}} = \PMod(\Sigma)/[\PMod(\Sigma),
\PMod(\Sigma)]\) of \(\PMod(\Sigma)\) is cyclic. In fact,
\(\PMod(\Sigma)^{\operatorname{ab}}\) vanishes when \(g \ge 3\).

\begin{theorem}[Powell, \cite{powell-78}]\label{thm:ab-vanishes}
  Let \(\Sigma\) be a surface of genus \(g \ge 3\). Then \(\PMod(\Sigma)\) is a
  perfect group.
\end{theorem}

The groups \(\PMod(\Sigma)^{\operatorname{ab}}\) are also known in the
low-genus cases. See \cite{korkmaz-02, farb-margalit}.

In a complementary direction, we can also consider the absence of relations
between two Dehn twists.

\begin{enumerate}[start=4]
  \item \textbf{Free subgroups.} Given \(a, b \subset \Sigma\) with \(|a
    \pitchfork b| \ge 2\), \(T_a\) and \(T_b\) generate a rank-\(2\) free group
    in \(\PMod(\Sigma)\) \cite[Theorem~3.14]{farb-margalit}.
\end{enumerate}

It turns out free subgroups are ubiquitous in \(\PMod(\Sigma)\). This is
because of the so-called \emph{Tits-alternative}: a subgroup \(G \le
\PMod(\Sigma)\) is either virtually Abelian or it contains free groups
\cite[Theorem~A]{mccarthy-85}. In particular, up to taking powers, the subgroup
generated by two mapping classes is almost always a free group.

\begin{theorem}[Theorem~B \cite{mccarthy-85}]\label{thm:tits-alternative-weak}
  Let \(\Sigma\) be a surface of genus \(g \ge 2\) and \(f, g \in
  \PMod(\Sigma)\). Then we can find \(n, m \ge 1\) such that \(f^n, g^m\)
  either commute or generate a free group.
\end{theorem}

In \textsection\ref{sec:matrix-graphs} we will make use of the following result
in our proof of Theorem~\ref{thm:bigger-bounds}.

\begin{proposition}\label{thm:normal-subgrps-contain-free}
  Let \(g \ge 2\) and \(N, K \triangleleft \Mod(\Sigma_g^b)\) be non-central
  normal subgroups. Then \(N \cap K\) contains a copy of the rank-\(2\) free
  group.
\end{proposition}

The idea is to show that any non-central \(N \triangleleft \Mod(\Sigma_g^b)\)
contains \emph{pseudo-Anosov maps}\footnote{See
\cite[\textsection13.2.3]{farb-margalit} for a definition.}, which is known
since the mid 1980's -- see \cite[Lemma~2.5]{long-86}. This is also a
consequence of the existence of all pseudo-Anosov normal subgroups of
\(\Mod(\Sigma_g^b)\), a fact first established by Dahmani--Guirardel--Osin
\cite[Theorem~2.31]{dahmani-17}.

\subsection{The Torelli group \& its Subgroups}
\label{sec:torelli-johnson}

Let \(\Sigma = \Sigma_{g, r}^b\) be a genus \(g\) surface with \(r\) marked
point \(P \subset \Sigma\) and \(b\) boundary components.

The natural action of \(\Diff(\Sigma, \partial \Sigma)_P\) on \(\Sigma\)
induces a \(\ZZ\)-linear action of \(\Mod(\Sigma)\) on the first homology group
\(H_1(\Sigma \setminus P; \ZZ)\). The Dehn twist \(T_a\) about \(a \subset
\Sigma\) acts by the operator
\begin{equation}\label{eq:homological-rep-formula}
  (T_a)_* x = x + \langle a, x \rangle a \in H_1(\Sigma \setminus P; \ZZ),
\end{equation}
where \(\langle \, , \rangle\) denotes the intersection pairing. Here we view
\(a \subset \Sigma\) as an element of \(H_1(\Sigma \setminus P; \ZZ)\) by
choosing an orientation of this curve. Notice, however, that
(\ref{eq:homological-rep-formula}) is independent of such a choice.

The \emph{Torelli subgroup} of \(\Sigma\), denoted \(\mathcal{I}(\Sigma) \le
\PMod(\Sigma)\), is the subgroup of mapping classes acting trivially on
\(H_1(\Sigma \setminus P; \ZZ)\). Some of its elements include the following
maps.

\begin{enumerate}
  \item \textbf{Bounding pair maps.} Given disjoint nonseparating curves \(a, b
    \subset \Sigma\) such that \(a \cup b\) bounds a closed subsurface
    \(\Sigma'\) with no marked points and boundary \(\partial \Sigma' = a \cup
    b\), the pair \((a, b)\) is called \emph{a bounding pair} and
    \(T_aT_b^{-1}\) is called \emph{a bounding pair map}. By
    (\ref{eq:homological-rep-formula}), \((T_a)_*\) only depends on the
    homology class of \(a\). Thus \(T_aT_b^{-1} \in \mathcal{I}(\Sigma)\).

  \item \textbf{Separating Dehn twists.} Given some separating \(a \subset
    \Sigma\), \(T_a\) is called a \emph{genus \(h\) separating Dehn twists} if
    \(a\) cuts \(\Sigma\) into subsurfaces of genus \(h\) and \(h'\) with \(h
    \le h'\). Since \(a\) is separating, \(\langle a, x \rangle = 0\) for all
    \(x \in H_1(\Sigma \setminus P; \ZZ)\). Hence \(T_a \in
    \mathcal{I}(\Sigma)\) by (\ref{eq:homological-rep-formula}).

  \item \textbf{Simple intersection maps.} Given \(a, b \subset \Sigma\) with
    geometric intersection number \(2\) and algebraic intersection number
    \(\langle a, b \rangle = 0\), the pair \((a, b)\) is called \emph{a simple
    intersection pair} and the commutator \([T_a, T_b] = T_a T_b T_a^{-1}
    T_b^{-1}\) is called \emph{a simple intersection map}.
    Since \(\langle a, b \rangle = 0\),
    \((T_a)_*\) and \((T_b)_*\) commute
    by
    (\ref{eq:homological-rep-formula}), so
    that \([T_a, T_b] \in \mathcal{I}(\Sigma)\).
\end{enumerate}

\begin{figure}
  \centering
  \begin{tikzpicture}[scale=1.25]
    \foreach \x/\i in {0/1, 1.5/2, 4/g} {
      \draw[Red, thick]                 (\x,0)+(-90:.35) arc[x radius=.6, y radius=.35, start angle=-90, end angle=270];
      \draw[Red, thick, -{>[scale=.5]}] (\x,0)+(-90:.35) arc[x radius=.6, y radius=.35, start angle=-90, end angle=270];
      \draw[Red, thick]                 (\x,0)+(-90:.35) node[below] {$b_{\i}$};
      \draw[blue, thick, line cap=round]         (\x,0)+(0,1)          arc[x radius=.1, y radius=.485, start angle=90,  end angle=-90];
      \draw[blue, thick, -{>[scale=.5]}]         (\x,0)+(0,{1-2*.485}) arc[x radius=.1, y radius=.485, start angle=-90, end angle=0];
      \draw[blue, thick, line cap=round, dotted] (\x,0)+(0,1)          arc[x radius=.1, y radius=.485, start angle=90,  end angle=270];
      \draw[blue]                                (\x,0)+(.1,{1-.485})  node[right] {$a_{\i}$};
      \torushole{\x}{0}
    }
    \foreach \i in {0,1,2} { \filldraw (2.5, 0)+({\i*.25}, 0) circle(0.5pt); }
    \draw[thick] (-.5,1) -- (4.5, 1) arc (90:-90:1) -- (-.5,-1) arc (-90:-270:1);
  \end{tikzpicture}
  \caption{A basis for the homology of $\Sigma_g$.}
  \label{fig:generators-of-homology}
\end{figure}
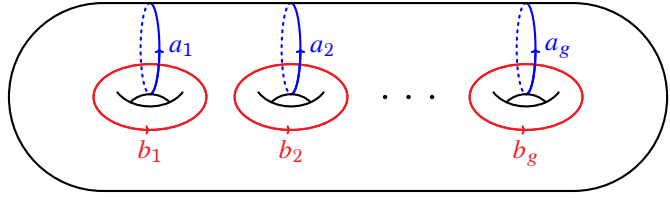

We can view \(\Sigma\) as a subsurface of the closed unmarked genus \(g\)
surface \(\Sigma_g\) by capping all boundary components with disks. Hence we
can also consider the action of \(\PMod(\Sigma)\) on \(H_\ZZ = H_1(\Sigma_g;
\ZZ)\), the first homology of the \emph{closed surface} of the same genus. This
is a rank-\(2g\) free Abelian group, freely generated by the homology class of
the oriented curves \(a_1, \ldots, a_g, b_1, \ldots, b_g\) in
Figure~\ref{fig:generators-of-homology}.

We thus obtain a linear representation \(\PMod(\Sigma) \to \GL_{2g}(\ZZ)\). Its
image is exactly the discrete symplectic group \(\Sp_{2g}(\ZZ)\). Indeed, the
action of \(\PMod(\Sigma)\) preserves the intersection pairing \(\langle \, ,
\rangle\), which corresponds to the standard symplectic form in \(H_1(\Sigma_g;
\RR) \cong \RR^{2g}\).

This is called the \emph{symplectic representation} \(\Psi : \PMod(\Sigma) \to
\Sp_{2g}(\ZZ)\). More concretely, \(\Psi\) is given by
\begin{equation}\label{eq:symplectic-rep-formula}
  \Psi(T_a) x = x + \langle a, x \rangle a \in H_\ZZ,
\end{equation}
not to be confused with the operators \((T_a)_*\) from
(\ref{eq:homological-rep-formula}). When \(\Sigma = \Sigma_g\),
\(\Psi(T_{a_i})\) and \(\Psi(T_{b_i})\) are given by the transvections
\begin{align*}
  \Psi(T_{a_i}) &=
  \left(
  \begin{array}{c|c|c}
    1 & 0                                          & 0 \\ \hline
    0 & \begin{matrix} 1 & 1 \\ 0 & 1 \end{matrix} & 0 \\ \hline
    0 & 0                                          & 1
  \end{array}
  \right) &
  \Psi(T_{b_i}) &=
  \left(
  \begin{array}{c|c|c}
    1 & 0                                           & 0 \\ \hline
    0 & \begin{matrix} 1 & 0 \\ -1 & 1 \end{matrix} & 0 \\ \hline
    0 & 0                                           & 1
  \end{array}
  \right)
\end{align*}
in the basis from Figure~\ref{fig:generators-of-homology}, where the top-left
blocks are $2(i - 1) \times 2(i-1)$.

The symplectic kernel \(\ker \Psi\) contains the Torelli subgroup
\(\mathcal{I}(\Sigma)\). This two subgroups coincide when \(\Sigma\) has no
marked points. Improving results of Birman \cite{birman-71} and Powell
\cite{powell-78}, Johnson \cite{johnson-83} found an explicit generating set
for the Torelli subgroup.

\begin{theorem}[Johnson, \cite{johnson-83}]
  Let \(\Sigma = \Sigma_g\) or \(\Sigma = \Sigma_g^1\) be an unmarked surface
  of genus \(g \ge 3\) with at most one boundary component. The Torelli
  subgroup \(\mathcal{I}(\Sigma)\) is generated by finitely many bounding pair
  maps.
\end{theorem}

Now assume \(\Sigma = \Sigma_g^1\) is an unmarked surface with one boundary
component. By choosing a base point \(* \in \partial \Sigma_g^1\), the mapping
class group \(\Mod(\Sigma_g^1)\) naturally acts on \(\Gamma = \pi_1(\Sigma_g^1,
*)\) by group automorphisms. The Torelli subgroup can then be seen as the
subgroup of mapping classes acting trivially on the Abelianization
\(\Gamma/[\Gamma, \Gamma] = H_1(\Sigma_g^1; \ZZ)\).

By considering the remaining terms \(\Gamma_k = [\Gamma, \Gamma_{k-1}]\) of the
lower central series of \(\Gamma = \Gamma_1\), Johnson introduced a filtration
\[
  \mathcal{I}(\Sigma_g^1) = \mathcal{I}^1(\Sigma_g^1)
  \triangleright \mathcal{I}^2(\Sigma_g^1)
  \triangleright \cdots
  \triangleright \mathcal{I}^k(\Sigma_g^1)
  \triangleright \cdots,
\]
known as the \emph{Johnson filtration} of \(\mathcal{I}(\Sigma_g^1)\). Here
\(\mathcal{I}^k(\Sigma_g^1)\) denotes the (normal) subgroup of mapping classes
acting trivially on the characteristic quotient \(\Gamma/\Gamma_k\).

The subgroups \(\mathcal{I}^k(\Sigma_g^1)\) are nontrivial for arbitrarily
large \(k\), and each \(\mathcal{I}^k(\Sigma_g^1)\) contains the
\(k^{\text{th}}\) term \(\mathcal{I}(\Sigma_g^1)_k\) of the lower central
series of \(\mathcal{I}(\Sigma_g^1) = \mathcal{I}(\Sigma_g^1)_1\). In
particular, \(\mathcal{I}^2(\Sigma_g^1) \ge [\mathcal{I}(\Sigma_g^1),
\mathcal{I}(\Sigma_g^1)]\) and the quotient
\(\mathcal{I}(\Sigma_g^1)/\mathcal{I}^2(\Sigma_g^1)\) is Abelian. Indeed,
Johnson essentially showed \(\mathcal{I}(\Sigma_g^1)/\mathcal{I}^2(\Sigma_g^1)
\cong {\bigwedge}^3 H_\ZZ\).

The subgroup \(\mathcal{K}(\Sigma_g^1) = \mathcal{I}^2(\Sigma_g^1)\) is called
the \emph{Johnson subgroup} of \(\mathcal{I}(\Sigma_g^1)\), while the projection
\(\tau : \mathcal{K}(\Sigma_g^1) \twoheadrightarrow {\bigwedge}^3 H_\ZZ\) is
called the \emph{Johnson homomorphism}. We may also consider
\(\mathcal{K}(\Sigma_g) = \ker \tau\), where \(\tau : \mathcal{I}(\Sigma_g)
\twoheadrightarrow {\bigwedge}^3 H_\ZZ/H_\ZZ\) is given by
\[
  \begin{tikzcd}
    \mathcal{I}(\Sigma_g^1)      \rar[two heads]{\tau} \dar            &
    {\bigwedge}^3 H_\ZZ                                \dar[two heads] \\
    \mathcal{I}(\Sigma_g)        \rar[two heads]{\tau}                 &
    {\bigwedge}^3 H_\ZZ / H_\ZZ.
  \end{tikzcd}
\]

Here the inclusion \(H_\ZZ \hookrightarrow {\bigwedge}^3 H_\ZZ\) takes \(x \in
H_\ZZ\) to \(a_1 \wedge b_1 \wedge x + \cdots + a_g \wedge b_g \wedge x\) for
\(a_1, b_1, \ldots, a_g, b_g\) as in Figure~\ref{fig:generators-of-homology}.
In an abuse of notation, \(\mathcal{K}(\Sigma_g)\) and \(\tau :
\mathcal{I}(\Sigma_g) \twoheadrightarrow {\bigwedge}^3 H_\ZZ/H_\ZZ\) are also
called the \emph{Johnson subgroup} and the \emph{Johnson homomorphism},
respectively. Johnson showed that \(\mathcal{K}(\Sigma_g)\) may also be
characterized as follows.

\begin{theorem}[Johnson, \cite{johnson-85}]
  Let \(g \ge 3\). Then \(\mathcal{K}(\Sigma_g)\) is the subgroup generated by
  all genus \(1\) and genus \(2\) separating Dehn twists.
\end{theorem}

The subgroup \(\SIP(\Sigma) \le \mathcal{I}(\Sigma)\) generated by simple
intersection maps remains less well understood than \(\mathcal{I}(\Sigma)\) and
\(\mathcal{K}(\Sigma)\). This is a normal subgroup of \(\Mod(\Sigma)\): given a
simple intersection pair \((a, b)\) and \(f \in \Mod(\Sigma)\),
\[
  f [T_{a}, T_{b}] f^{-1}
  = [f T_{a} f^{-1}, f T_{b} f^{-1}]
  = [T_{f(a)}, T_{f(b)}] \in \SIP(\Sigma)
\]
by the conjugation relation.

When \(\Sigma = \Sigma_g^1\), the images of simple intersection maps under the
Johnson homomorphism were computed independently by Putman \cite{putman-15},
Church \cite{church-14} and Childers \cite[Main Result~1]{childers-12}.
Childers also computed the image of \(\SIP(\Sigma_g^1)\) under the so-called
\emph{Birman--Craggs--Johnson homomorphism} \cite[Main Result~4]{childers-12}.

We further restrict our attention to the subgroup \(\SIP_0(\Sigma) \le
\SIP(\Sigma)\) generated by simple intersection maps \([T_a, T_b]\) where
\(\Sigma \setminus (a \cup b)\) is connected. This is also a normal subgroup,
normally generated by \([T_a, T_b]\) for \emph{any} simple intersection pair
\((a, b)\) with \(\Sigma \setminus (a \cup b)\) connected. Indeed, given any
other choice of \((a', b')\) as above, we can find \(f \in \Mod(\Sigma)\) with
\(f(a) = a'\) and \(f(b) = b'\), implying that all generators of
\(\SIP_0(\Sigma)\) are conjugate in \(\Mod(\Sigma)\).

\subsection{The Work of Korkmaz}

Let \(\Sigma\) be a surface as above. We now describe some results due to
Korkmaz which are needed in the rest of the article. We begin by his
classification theorem.

The starting point of Korkmaz' classification program is the aforementioned
computation of the Abelianization of \(\PMod(\Sigma)\). Using the braid and
disjointness relations, Korkmaz \cite{korkmaz-23} showed that, when \(g \ge
2\), any linear representation \(\PMod(\Sigma) \to \GL_d(\CC)\) with \(d < 2g\)
factors through the Abelianization map \(\PMod(\Sigma) \twoheadrightarrow
\PMod(\Sigma)^{\operatorname{ab}}\).

Since \(\PMod(\Sigma)^{\operatorname{ab}}\) vanishes for \(g \ge 3\)
(Theorem~\ref{thm:ab-vanishes}), any such representation must be trivial.
Korkmaz furthermore used the same relations to show that all nontrivial
\(2g\)-dimensional representations are conjugate to the symplectic
representation \(\Psi\).

\begin{theorem}[Korkmaz, Theorems~1 \& 2 \cite{korkmaz-23}]\label{thm:korkmaz-main-thms}
  Let \(\Sigma\) be a surface of genus \(g \ge 3\) and \(\rho : \PMod(\Sigma)
  \to \GL_d(\CC)\). If \(d \le 2g\) then \(\rho\) is either trivial or
  conjugate to \(\Psi : \PMod(\Sigma) \to \Sp_{2g}(\ZZ)\).
\end{theorem}

Combining this last result with Theorem~\ref{thm:ab-vanishes}, Korkmaz also
established the following triviality criterion.

\begin{lemma}[Flag triviality criterion, Lemma~7.1 \cite{korkmaz-23}]\label{thm:korkmaz-triviality-criteria}
  Let \(\Sigma\) be a surface of genus \(g \ge 3\) and \(\rho : \PMod(\Sigma)
  \to \GL_d(\CC)\). Suppose there exists a \(\PMod(\Sigma)\)-invariant flag
  \[
    0 \le W_1 \le W_2 \le \cdots \le W_k = \CC^d
  \]
  with \(\dim W_k/W_{k+1} < 2g\). Then \(\rho\) is trivial.
\end{lemma}

As it turns out, invariant flags are pervasive. This is because of the
following principle. Given \(\rho : \PMod(\Sigma) \to \GL_d(\CC)\) and \(a
\subset \Sigma\), denote \(L_a = \rho(T_a)\). If \(a, b \subset \Sigma\) are
disjoint then \(L_a\) and \(L_b\) commute, so that \(L_b\) preserves the
eigenspaces of \(L_a\). Combining this observation with
Theorem~\ref{thm:mcg-generators} we obtain the following.

\begin{lemma}[Korkmaz, Lemma~4.1 \cite{korkmaz-23}]\label{thm:eigenspace-is-invariant}
  Let \(\Sigma' \subset \Sigma\) be closed subsurface and \(\rho :
  \PMod(\Sigma) \to \GL_d(\CC)\). Take \(a \subset \Sigma \setminus \Sigma'\).
  Then \(E_{\lambda, k}^a = \ker (L_a - \lambda)^k\) is a
  \(\PMod(\Sigma')\)-invariant subspace of \(\CC^d\). In particular, the flag
  \[
    0
    \le E_{\lambda, 1}^a
    \le E_{\lambda, 2}^a
    \le \cdots
    \le E_{\lambda, d}^a
    \le \CC^d
  \]
  is \(\PMod(\Sigma')\)-invariant.
\end{lemma}

\subsection{Twisted Cohomology}
\label{sec:twisted-cohomology}

Let \(\Sigma = \Sigma_{g, r}^b\) be the genus \(g\) compact surface with \(b\)
boundary components and \(r\) marked points. In this subsection we review some
results on the cohomology of \(\PMod(\Sigma)\). These will be used in
\textsection\ref{sec:eigenvals} to handle the \(d = 4g - 4\) case of
Theorem~\ref{thm:bound-on-faithful-mcg}. We refer the reader to \cite{brown}
for a comprehensive account of the theory of group cohomology.

Denote by \(\ZZ[\PMod(\Sigma)] = \bigoplus_{f \in \PMod(\Sigma)} \ZZ f\) the
\emph{group ring} of \(\PMod(\Sigma)\): the ring of (formal) integral
combinations of elements in \(\PMod(\Sigma)\), where multiplication is given by
the product in \(\PMod(\Sigma)\). Given a representation \(\rho : \PMod(\Sigma)
\to \GL_d(R)\), we may view \(R^d\) as a \(\ZZ[\PMod(\Sigma)]\)-module where
\(f \in \PMod(\Sigma)\) acts by \(\rho(f)\). For example, the \emph{trivial}
\(\ZZ[\PMod(\Sigma)]\)-module is the module \(\ZZ\) corresponding to the
trivial homomorphism \(\PMod(\Sigma) \to \GL_1(\ZZ)\).

Recall that, given a \(\ZZ[\PMod(\Sigma)]\)-module \(M\), a map \(c :
\PMod(\Sigma) \to M\) is called an \emph{\(M\)-valued crossed homomorphism} if
\(c(f g) = c(f) + f \cdot c(g)\) for all \(f, g \in \PMod(\Sigma)\). The
collection of all such maps forms an Abelian group. A crossed homomorphism
\(c\) is called \emph{principal} if there is \(m \in M\) such that \(c(f) = m -
f \cdot m\) for all \(f\).

The \emph{first group cohomology group of \(\PMod(\Sigma)\) with coefficients
in \(M\)}, denoted \(H^1(\PMod(\Sigma); M)\), is the quotient of the group of
crossed homomorphisms by the subgroup of principal crossed homomorphisms. Its
elements are in one-to-one correspondence with isomorphism classes of
extensions of the trivial \(\ZZ[\PMod(\Sigma)]\)-module by \(M\), i.e. short
exact sequences of the form
\[
  \begin{tikzcd}
    0 \rar & M \rar & E \rar & \ZZ \rar & 0.
  \end{tikzcd}
\]

In a similar way, one may define \emph{higher cohomology groups}
\(H^k(\PMod(\Sigma); M)\) for \(k \ge 0\).

As a consequence of Theorem~\ref{thm:ab-vanishes},
we obtain the following computation.

\begin{lemma}\label{thm:first-cohomology-vanishes}
  Let \(\Sigma\) be a surface of genus \(g \ge 3\). Then \(H^1(\PMod(\Sigma);
  \ZZ) = 0\).
\end{lemma}

We also consider the cohomology with coefficients in the
\(\ZZ[\PMod(\Sigma)]\)-modules \(H_\ZZ = H_1(\Sigma_g; \ZZ)\) and \(H_\CC =
H_1(\Sigma_g; \CC)\) corresponding to the symplectic representation \(\Psi :
\PMod(\Sigma) \to \Sp_{2g}(\ZZ)\). We pay special attention to the surfaces
\(\Sigma_g^1\) and \(\Sigma_{g, 1}\). Such groups were first computed by Morita
\cite{morita-89}.

\begin{theorem}[Morita, Proposition~6.4 \cite{morita-89}]
  If \(g \ge 2\) then \(H^1(\Mod(\Sigma_g^1); H_\ZZ) \cong H^1(\Mod(\Sigma_{g,
  1}); H_\ZZ) \cong \ZZ\).
\end{theorem}

A standard argument with the universal coefficient theorem shows
\(H^1(\Mod(\Sigma_g^1); H_\CC) \cong \CC\) for \(g \ge 2\). This implies that
all nontrivial extensions of the \(H_\CC\) by \(\CC\) are isomorphic as
\(\ZZ[\Mod(\Sigma_g^1)]\)-modules. Kasahara \cite{kasahara-23} showed that, up
to dualizing, any nontrivial \((2g + 1)\)-dimensional \(\ZZ[\PMod(\Sigma_{g,
r}^b)]\)-module is isomorphic to one such extension. Here we make use of a
homological lemma of his.

\begin{lemma}[Kasahara, Theorem~4.2 \cite{kasahara-23}]\label{thm:crossed-hom-characterization}
  Let \(c : \Mod(\Sigma_g^1) \to H_\CC\) be a crossed homomorphism. Given \(a
  \subset \Sigma_g^1\) nonseparating, \(c(T_a) = \lambda \cdot a\) for some
  \(\lambda \in \CC\).
\end{lemma}

Morita's result was later generalized by Kawazumi \cite{kawazumi-08}, who
computed the (stable) higher cohomology groups of \(\Mod(\Sigma_g^1)\) with
coefficients in tensor powers of \(H_\ZZ\) in terms of the so called
\emph{twisted Miller--Mumford--Morita classes} \(\hat\kappa_P \in
H^*(\Mod(\Sigma_g^1); H_\ZZ^{\otimes n})\) associated to \emph{weighted
partitions} \(P\) of the set \(\{1, \ldots, n\}\). See \cite{kawazumi-08} for
definitions.

\begin{theorem}[Kawazumi, Theorem~1.B \cite{kawazumi-08}]\label{thm:tensor-prod-cohomology}
  For \(k \le \sfrac{g}{2} - n\),
  \[
    H^k(\Mod(\Sigma_g^1); H_\ZZ^{\otimes n})
    = \bigoplus_{\substack{P \in \mathcal{P}_n \\ \ell + \deg \hat\kappa_P = k}}
      H^\ell(\Mod(\Sigma_g^1); \ZZ) \smile \hat\kappa_P.
  \]
\end{theorem}

Here \(\mathcal{P}_n\) denotes the set of weighted partitions of \(\{1, \ldots,
n\}\) and
\[
  H^\ell(\Mod(\Sigma_g^1); \ZZ) \smile \hat\kappa_P
  = \{\xi \smile \hat\kappa_P : \xi \in H^\ell(\Mod(\Sigma_g^1); \ZZ)\},
\]
where \(\xi \smile \hat\kappa_P \in H^k(\Mod(\Sigma_g^1); H_\ZZ^{\otimes n})\)
is the \emph{cup product} of \(\xi\) with \(\hat\kappa_P\).

When \(n = 2\), the twisted Miller--Mumford--Morita classes take two forms: the
classes \(\alpha_i \in H^{2i}(\Mod(\Sigma_g^1); H_\ZZ^{\otimes 2})\) where \(i
\ge 0\), and the classes \(\beta_{ij} \in H^{2i + 2j - 2}(\Mod(\Sigma_g^1);
H_\ZZ^{\otimes 2})\) where \(i, j \ge 1\), corresponding to partitions of
\(\{1, 2\}\) into one and two subsets, respectively. We thus obtain the
following computation.

\begin{corollary}\label{thm:coh-tensor-square-vanishes}
  For \(g \ge 6\), \(H^1(\Mod(\Sigma_g^1); H_\ZZ^{\otimes 2}) = 0\).
\end{corollary}

\begin{proof}
  Since \(1 \le \sfrac{g}{2} - 2\), Theorem~\ref{thm:tensor-prod-cohomology}
  says \(H^1(\Mod(\Sigma_g^1); H_\ZZ^{\otimes 2}) = H^1(\Mod(\Sigma_g^1); \ZZ)
  \smile \alpha_0\). Now recall from Lemma~\ref{thm:first-cohomology-vanishes}
  that \(H^1(\Mod(\Sigma_g^1); \ZZ) = 0\) already for \(g \ge 3\).
\end{proof}

\section{Commutation Relations in the Mapping Class Group}
\label{sec:matrix-graphs}

Recall that the \emph{commuting graph} of a group \(G\) is the graph
\(\Gamma(G)\) whose vertices are elements of \(G\), where \(g, h \in G\) are
joined by an edge if and only if they commute. First defined by Harvey, the
\emph{curve graph} of a surface \(\Sigma\) is the graph \(\mathcal{C}(\Sigma)\)
whose vertices are homotopy classes of essential simple closed curves in
\(\Sigma\), where \(a, b \subset \Sigma\) are joined by an edge if and only if
we can find disjoint representatives.

The latter is a Gromov-hyperbolic graph on which \(\PMod(\Sigma)\) acts by
isometries \cite{masur-99}. Given the disjointness relations from
\textsection\ref{sec:dehn-twists}, \(\mathcal{C}(\Sigma)\) and
\(\Gamma(\PMod(\Sigma))\) are related by means of the embedding
\(\mathcal{C}(\Sigma) \hookrightarrow \Gamma(\PMod(\Sigma))\) taking \(a
\subset \Sigma\) to its Dehn twist \(T_a\).

The starting point of the present article was to consider the disjointness
relations in \(\PMod(\Sigma)\) induced by a family curves \(a_1, \ldots, a_n,
b_1, \ldots, b_n \subset \Sigma\) satisfying
\begin{align}\label{eq:intersection-rels}
  |a_i \pitchfork a_j| & = 0   \;\; \forall i \ne j &
  |b_i \pitchfork b_j| & = 0   \;\; \forall i \ne j &
  |a_i \pitchfork b_j| & = 0   \;\; \forall i \ne j &
  |a_i \pitchfork b_i| & \ge 1 \;\; \forall i,
\end{align}
where \(|a \pitchfork b|\) denotes the geometric intersection number between
\(a\) and \(b\).

\begin{figure}
  \centering
  \begin{tikzpicture}
    \foreach \i in {1, 2, 3, 4} {
      \filldraw (0,-\i) circle(2pt);
      \filldraw (2,-\i) circle(2pt);
    }
    \filldraw (0,-1) node[anchor=south east] {$a_1$}
              (2,-1) node[anchor=south west] {$b_1$}
              (0,-2) node[anchor=south east] {$a_2$}
              (2,-2) node[anchor=south west] {$b_2$};
    \filldraw (0,-3) node[anchor=north east] {$a_{n-1}$}
              (2,-3) node[anchor=north west] {$b_{n-1}$}
              (0,-4) node[anchor=north east] {$a_n$}
              (2,-4) node[anchor=north west] {$b_n$};
    \foreach \i/\j in {1/2,1/3,1/4,2/1,2/3,2/4,3/1,3/2,3/4,4/1,4/2,4/3}
      \draw[thick] (0,-\i) -- (2,-\j);
    \draw[thick]         (0, -1) -- (0, -2);
    \draw[thick, dotted] (0, -2) -- (0, -3);
    \draw[thick]         (0, -3) -- (0, -4);
    \draw[thick]         (2, -1) -- (2, -2);
    \draw[thick, dotted] (2, -2) -- (2, -3);
    \draw[thick]         (2, -3) -- (2, -4);
    \draw[thick] (0,-1) arc (90:270:1);
    \draw[thick] (0,-2) arc (90:270:1);
    \draw[thick] (0,-1) arc[x radius=1.25, y radius=1.5, start angle=90, end angle=270];
    \draw[thick] (2,-1) arc (90:-90:1);
    \draw[thick] (2,-2) arc (90:-90:1);
    \draw[thick] (2,-1) arc[x radius=1.25, y radius=1.5, start angle=90, end angle=-90];
  \end{tikzpicture}
  \caption{The graph $\Delta_n$: the vertices $a_1, \ldots, a_n$ and $b_1, \ldots, b_n$ form two disjoint $n$-cliques, while the vertices $a_i$ and $b_j$ are connected by an edge if and only if $i \ne j$.}
  \label{fig:delta-graph}
\end{figure}
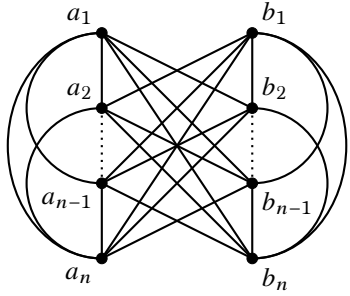

These correspond to copies of the graph \(\Delta_n\) from
Figure~\ref{fig:delta-graph} inside \(\mathcal{C}(\Sigma)\). In terms of
\(\Gamma(\PMod(\Sigma))\), such a family of curves translates to the relations
\begin{align}\label{eq:comm-rels-mcg}
  T_{a_i} T_{a_j} & = T_{a_j} T_{a_i} \text{ for all } i, j &
  T_{b_i} T_{b_j} & = T_{b_j} T_{b_i} \text{ for all } i, j &
  T_{a_i} T_{b_j} & = T_{b_j} T_{a_i} \iff i \ne j.
\end{align}
These relations can be further refined as follows.

Assume \(a_1, \ldots, a_n, b_1, \ldots, b_n \subset \Sigma\) is a family as
above, with the first \(k\) pairs \((a_i, b_i)\) intersecting exactly once and
the remaining \(\ell = n - k\) pairs intersecting at \(\ge 2\) points. By the
braid relation, \(T_{a_i}\) and \(T_{b_i}\) generate a copy of the braid group
\(B_3\) inside \(\PMod(\Sigma)\) for \(i \le k\). Likewise, \(T_{a_i}\) and
\(T_{b_i}\) generate a rank-\(2\) free group \(F_2\) for \(i \ge k + 1\).

Combined with (\ref{eq:comm-rels-mcg}), these relations then imply there is a a
well defined homomorphism from
\[
  G_{k\ell} =
  \underbrace{B_3 \times \cdots \times B_3}_{k \text{ times}}
  \times
  \underbrace{F_2 \times \cdots \times F_2}_{\ell \text{ times}}
\]
onto the subgroup generated by \(T_{a_1}, \ldots, T_{a_n}, T_{b_1}, \ldots,
T_{b_n}\). What is more, the projection of each direct factor of \(G_{k\ell}\)
onto its image in \(\PMod(\Sigma)\) is an isomorphism.

More generally, one can produce quotients of \(G_{k\ell}\) inside a subgroup
\(G \le \PMod(\Sigma)\) by considering the subgroups of \(G\) consisting of
mapping classes supported on disjoint closed subsurfaces \(S_{1}, \ldots, S_{n}
\subset \Sigma\), each one corresponding to one of the factors of
\(G_{k\ell}\). It is thus natural to expect the dimension of a faithful
representation \(G \hookrightarrow \GL_d(\CC)\) to be related to the
minimal dimension of a faithful representation of \(G_{k\ell}\).

The latter was recently estimated by Kionke--Schesler
\cite{kionke-schesler-23}.

\begin{proposition}[Kionke--Schesler, Theorem~3
  \cite{kionke-schesler-23}]\label{thm:non-solvable-prod-linear-bound}
  Suppose \(H_1, \ldots, H_n\) are non-solvable groups and \(\rho : H_1 \times
  \cdots \times H_n \hookrightarrow \GL_d(\CC)\) is a faithful linear
  representation. Then \(d \ge 2n\).
\end{proposition}

Their proof is short and elementary, making clever use of well known facts
about the representation theory of direct products. Here we make use of a
slightly more general version of their statement, although our proof is really
an adaptation of their argument.

\begin{proposition}\label{thm:free-group-prod-linear-bound}
  Suppose \(H_1, \ldots, H_n\) are non-solvable groups and \(\pi : H_1 \times
  \cdots \times H_n \twoheadrightarrow H\) is a surjective group homomorphism
  such that \(\pi\!\restriction_{H_i} : H_i \to H\) is injective for all \(i\).
  Let \(\rho : H \hookrightarrow \GL_d(\CC)\) be a faithful linear
  representation. Then \(d \ge 2n\).
\end{proposition}

\begin{proof}
  Consider \(\rho \circ \pi : H_1 \times \cdots \times H_n \to \GL_d(\CC)\) and
  take a maximal \((H_1 \times \cdots \times H_n)\)-invariant flag
  \[
    0 = W_0 \le W_1 \le \cdots \le W_{p+1} = \CC^d,
  \]
  so that the action \(\rho_i : H_1 \times \cdots \times H_n \to
  \GL(W_{i+1}/W_i)\) of \(H_1 \times \cdots \times H_n\) on each successive
  quotient is irreducible. In a basis adapted to this flag,
  \begin{equation}\label{eq:kionke-schesler-basis}
    \rho(\pi(h)) =
    \begin{pmatrix}
      \rho_1(h) &         * & \cdots &         * \\
              0 & \rho_2(h) & \cdots &         * \\
         \vdots &    \vdots & \ddots &    \vdots \\
              0 &         0 & \cdots & \rho_p(h)
    \end{pmatrix}
  \end{equation}
  for all \(h \in H_1 \times \cdots \times H_n\).

  Now each \(W_{i+1}/W_i\) may be decomposed as a tensor product
  \(W_{i+1}/W_i = W_{i,1} \otimes \cdots \otimes W_{i,n}\), where
  \(\rho_{i j} : H_j \to \GL(W_{ij})\) is an irreducible representation and
  \(\rho_i(h_1, \ldots, h_n) = \rho_{i,1}(h_1) \otimes \cdots \otimes
  \rho_{i,n}(h_n)\) -- see, for example, \cite[Proposition~2.3.23]{kowalski}.
  In this setting,
  \[
    \dim W_{i+1}/W_i
    = (\dim W_{i, 1}) \cdots (\dim W_{i, n})
    \ge 2^{\# \{j : \dim W_{ij} \ge 2\}}
    \ge 2 \cdot \# \{j : \dim W_{ij} \ge 2\},
  \]
  for \(2^k \ge 2k\) for all integers \(k \ge 0\).
  In particular, \(d = \dim W_1/W_0 + \cdots + \dim W_{p+1}/W_p \ge 2 \cdot \#
  \{(i, j) : \dim W_{ij} \ge 2 \}\).

  If \(\dim W_{i j} = 1\) then \(\rho_{ij}(h_j)\) is a scalar operator for all
  \(h_j \in H_j\), so that \(\rho_i(1, \ldots, h_j, \ldots, 1) = 1 \otimes
  \cdots \otimes \rho_{ij}(h_j) \otimes \cdots \otimes 1\) is also a scalar
  operator. Assume we can find \(j \le n\) with \(\dim W_{i j} = 1\) for all
  \(i\). Then the matrix \(\rho(\pi(1, \ldots, h_j, \ldots, 1))\) is upper
  triangular with respect to the basis from (\ref{eq:kionke-schesler-basis})
  for all \(h_j \in H_j\). Since the group of upper triangular matrices is
  solvable, it follows \(\pi(H_i^{(k)}) \le \ker \rho\) for some \(k\).

  But \(H_i\) is non-solvable and so \(\pi(H_i^{(k)}) \cong H_i^{(k)} \ne 1\),
  contradicting the assumption that \(\rho\) is faithful. This means that, for
  each \(j \le n\) we can find \(i\) such that \(\dim W_{i j} \ge 2\). Hence
  \(d \ge 2 \cdot \# \{(i, j) : \dim W_{ij} \ge 2 \} \ge 2n\), as desired.
\end{proof}

Our proofs of Theorem~\ref{thm:bound-on-faithful-johnson},
Theorem~\ref{thm:bound-on-faithful-johnson-filtration} and
Theorem~\ref{thm:bound-on-faithful-braid} are direct applications of
Proposition~\ref{thm:free-group-prod-linear-bound}.

\begin{corollary}[Theorem~\ref{thm:bound-on-faithful-johnson}]\label{thm:small-dim-rep-kills-separating-twists}
  Let \(g \ge 2\) and suppose \(\rho : \mathcal{K}(\Sigma_g) \hookrightarrow
  \GL_d(\CC)\) is a faithful representation of the Johnson kernel. Then \(d \ge
  2g - 2\).
\end{corollary}

\begin{proof}
  Consider the curves \(a_1, \ldots, a_{g-1}, c_1, \ldots, c_{g-1} \subset
  \Sigma_g\) as in Figure~\ref{fig:johnson-curves-def-closed} and take
  \(b_i = T_{c_i}(a_i)\). For each \(i\), the curves \(a_i\) and \(c_i\)
  intersect twice. The geometric intersection number of \(a_i\) and \(b_i\) is
  thus \(4 = 2^2\). It is also clear all other pairs of curves in the above
  family are disjoint.

  \begin{figure}
    \centering
    \begin{subfigure}[b]{\textwidth}
      \centering
      \begin{tikzpicture}[scale=1.25]
        \foreach \x/\i in {0/1, 3.75/g-1} {
          \draw[Red, thick, line cap=round]
            (\x, 0)+({1.75/2},.225) arc[start angle=90, end angle=0,   y radius=0.25, x radius=0.7]
            (\x, 0)+({1.75/2},.225) arc[start angle=90, end angle=180, y radius=0.25, x radius=0.7]
            node[above=.5em]{$c_{\i}$};
          \draw[Red, thick, dotted, line cap=round]
            (\x, 0)+({1.75/2},-.225) arc[start angle=-90, end angle=-20,    y radius=0.25, x radius=0.7]
            (\x, 0)+({1.75/2},-.225) arc[start angle=-90, end angle=-160, y radius=0.25, x radius=0.7];
          \draw[blue, thick, line cap=round]
            (\x, 0)+({1.75/2},1)
              arc[start angle=90, end angle=-90, y radius=1, x radius=0.2]
              node[below]{$a_{\i}$};
          \draw[blue, thick, dotted, line cap=round]
            (\x,0)+({1.75/2},1)
            arc[start angle=90, end angle=270, y radius=1, x radius=0.2];
        }
        \draw[thick] (-.5,1) -- (6, 1) arc (90:-90:1) -- (-.5,-1) arc (270:90:1);
        \torushole{0}{0}
        \torushole{1.75}{0}
        \torushole{3.75}{0}
        \torushole{5.5}{0}
        \foreach \i in {0,1,2} { \filldraw (2.5,0)+({\i*.25}, 0) circle(0.4pt); }
      \end{tikzpicture}
      \caption{$\Sigma = \Sigma_g$}
      \label{fig:johnson-curves-def-closed}
    \end{subfigure}
    \\[3ex]
    \begin{subfigure}[b]{\textwidth}
      \centering
      \begin{tikzpicture}[scale=1.25]
        \foreach \x/\i in {0/1, 3.75/g-1} {
          \draw[Red, thick, line cap=round]
            (\x, 0)+({1.75/2},.225) arc[start angle=90, end angle=0,   y radius=0.25, x radius=0.7]
            (\x, 0)+({1.75/2},.225) arc[start angle=90, end angle=180, y radius=0.25, x radius=0.7]
            node[above=.5em]{$c_{\i}$};
          \draw[Red, thick, dotted, line cap=round]
            (\x, 0)+({1.75/2},-.225) arc[start angle=-90, end angle=-20,    y radius=0.25, x radius=0.7]
            (\x, 0)+({1.75/2},-.225) arc[start angle=-90, end angle=-160, y radius=0.25, x radius=0.7];
          \draw[blue, thick, line cap=round]
            (\x, 0)+({1.75/2},1)
              arc[start angle=90, end angle=-90, y radius=1, x radius=0.2]
              node[below]{$a_{\i}$};
          \draw[blue, thick, dotted, line cap=round]
            (\x,0)+({1.75/2},1)
            arc[start angle=90, end angle=270, y radius=1, x radius=0.2];
        }
        \draw[thick] (-1.25,1) -- (6, 1) arc (90:-90:1) -- (-1.25,-1);
        \draw[thick, line cap=round] (-1.25,1) arc[start angle=90, end angle=-270, y radius=1, x radius=.2];
        \torushole{0}{0}
        \torushole{1.75}{0}
        \torushole{3.75}{0}
        \torushole{5.5}{0}
        \foreach \i in {0,1,2} { \filldraw (2.5,0)+({\i*.25}, 0) circle(0.4pt); }
      \end{tikzpicture}
      \caption{$\Sigma = \Sigma_g^1$}
      \label{fig:johnson-curves-def-boundary}
    \end{subfigure}
    \caption{The curves $a_1, \ldots, a_{g-1}, c_1, \ldots, c_{g-1} \subset \Sigma$.}
  \end{figure}
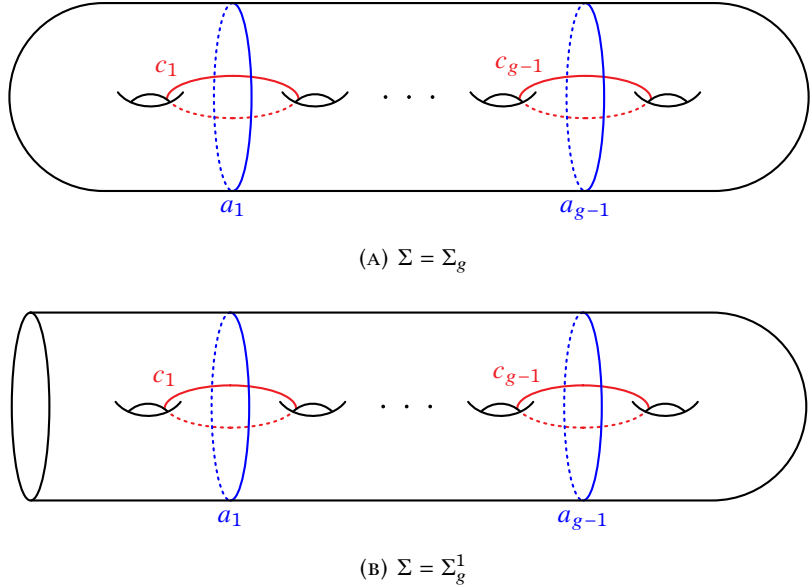

  The curves \(a_i\) and \(b_j\) are all separating. The discussion above then
  implies that the subgroup of \(\mathcal{K}(\Sigma_g)\) generated by
  \(T_{a_1}, \ldots, T_{a_{g-1}}, T_{b_1}, \ldots, T_{b_{g-1}}\) is a quotient
  of \(G_{0, g-1}\), the direct product of \(g-1\) copies of a rank-\(2\) free
  group \(F_2\). What is more, the projection of each \(F_2\)-factor of \(G_{0,
  g-1}\) onto its image in \(\mathcal{K}(\Sigma_g)\) is an isomorphism. The
  result thus follows from Proposition~\ref{thm:free-group-prod-linear-bound}.
\end{proof}

To prove Theorem~\ref{thm:bound-on-faithful-johnson-filtration}, we pass to the
derived subgroups \(\mathcal{I}(\Sigma_g^1)^{(k)} =
[\mathcal{I}(\Sigma_g^1)^{(k-1)}, \mathcal{I}(\Sigma_g^1)^{(k-1)}]\) of
\(\mathcal{I}(\Sigma_g^1) = \mathcal{I}(\Sigma_g^1)^{(1)}\).

\begin{corollary}[Theorem~\ref{thm:bound-on-faithful-johnson-filtration}]\label{thm:bound-on-faithful-johnson-filtration-real}
  Let \(\Sigma_g^1\) be the unmarked genus \(g\) surface with one boundary
  component. Suppose \(g \ge 2\). If \(\rho : \mathcal{I}(\Sigma_g^1)^{(k)}
  \hookrightarrow \GL_d(\CC)\) is faithful then \(d \ge 2g - 2\). In
  particular, if \(\rho : \mathcal{I}^k(\Sigma_g^1) \hookrightarrow
  \GL_d(\CC)\) is faithful then \(d \ge 2g - 2\).
\end{corollary}

\begin{proof}
  For the first claim, consider the curves \(a_1, \ldots, a_{g-1}, c_1, \ldots,
  c_{g-1} \subset \Sigma_g^1\) from
  Figure~\ref{fig:johnson-curves-def-boundary} and take \(b_i = T_{c_i}(a_i)\).
  As in the proof of Corollary~\ref{thm:small-dim-rep-kills-separating-twists},
  the subgroup of \(\mathcal{I}(\Sigma_g^1)\) generated by \(T_{a_1}, \ldots,
  T_{a_{g-1}}, T_{b_1}, \cdots, T_{b_{g-1}}\) is a quotient of \(G_{0,g-1} =
  F_2 \times \cdots \times F_2\).

  In particular, \(\mathcal{I}(\Sigma_g^1)^{(k)}\) contains a quotient of the
  \(k^{\text{th}}\) derived subgroup \(G_{0,g-1}^{(k)} = F_2^{(k)} \times
  \cdots \times F_2^{(k)}\). What is more, the projection of each
  \(F_2^{(k)}\)-factor onto its image in \(\mathcal{I}(\Sigma_g^1)^{(k)}\) is
  an isomorphism. Since \(F_2\) is non-solvable, so is \(F_2^{(k)}\). The result
  thus follows from Proposition~\ref{thm:free-group-prod-linear-bound}.

  For the second claim, it suffices to observe \(\mathcal{I}(\Sigma_g^1)^{(k)}
  \le \mathcal{I}(\Sigma_g^1)_k \le \mathcal{I}^k(\Sigma_g^1)\), where
  \(\mathcal{I}(\Sigma_g^1)_k = [\mathcal{I}(\Sigma_g^1),
  \mathcal{I}(\Sigma_g^1)_{k-1}]\) are the terms of the lower central series of
  \(\mathcal{I}(\Sigma_g^1) = \mathcal{I}(\Sigma_g^1)_1\).
\end{proof}

\begin{corollary}[Theorem~\ref{thm:bound-on-faithful-braid}]\label{thm:bound-on-faithful-braid-real}
  Suppose \(\rho : PB_n \hookrightarrow \GL_d(\CC)\) is faithful. If \(n\) is
  odd then \(d \ge n - 1\). If \(n\) is even then \(d \ge n - 2\).
\end{corollary}

\begin{proof}
  Recall \(PB_n = \PMod(\DD_n)\) is the pure mapping class group of
  a disk with \(n\) marked points. Given \(m < n\), the natural map \(PB_m \to
  PB_n\) is injective, so that we may regard \(PB_m\) as subgroup of \(PB_n\).

  In particular, if \(n\) is even, we can pass to the subgroup \(PB_{n-1} \le
  PB_n\). We may thus assume \(n = 2k + 1\) for some \(k \ge 0\). Furthermore,
  the result is clearly true for \(n \le 3\). We can thus assume \(k \ge 2\).

  \begin{figure}
    \centering
    \begin{tikzpicture}[scale=1.25]
      \draw[thick] (-2,0) arc[start angle=180, end angle=540, radius={(4+.75+2)/2}, y radius=1.7];
      \foreach \i in {0,1,2,3,4} { \filldraw (0,0)+({\i*.5}, 0) circle(1pt); }
      \foreach \i in {0,1,2} { \filldraw (2.5, 0)+({\i*.25}, 0) circle(0.5pt); }
      \foreach \i in {0,1} { \filldraw (3.5,0)+({\i*.5}, 0) circle(1pt); }
      \draw[thick, Red] (.25,0) ellipse (.4 and .2) +(-.4,0) node[left]{$b_1$};
      \draw[thick, Red]
        (-.75,0) arc[start angle=180, end angle=540, x radius={(1.5+.15+.75)/2}, y radius=.6]
        node[left]{$b_2$};
      \draw[thick, Red]
        (-1.25,0) arc[start angle=180, end angle=540, x radius={(3.5+.15+1.25)/2}, y radius=1.1]
        node[left]{$b_k$};
      \foreach \i/\l in {0/1,1/2,3/k} {
        \draw[thick, blue] (.75,0)+(\i,0) ellipse (.4 and .2) node[below=.75em]{$a_\l$};
      }
    \end{tikzpicture}
    \caption{The curves $a_1, \ldots, a_k, b_1, \ldots, b_k \subset \DD_n$.}
    \label{fig:braid-curves-def}
  \end{figure}

  In that case, consider \(a_1, \ldots, a_k, b_1, \ldots, b_k \subset \DD_n\)
  as in Figure~\ref{fig:braid-curves-def}. It follows from the discussion above
  that the subgroup generated by \(T_{a_1}, \ldots, T_{a_k}, T_{b_1}, \ldots,
  T_{b_k} \in \PMod(\DD_n)\) is a quotient of \(G_{0, k}\), the direct product
  of \(k\) copies of \(F_2\). What is more, the projection of each
  \(F_2\)-factor onto its image in \(\PMod(\Sigma_{0, n}^1)\) is an
  isomorphism. Proposition~\ref{thm:free-group-prod-linear-bound} then says \(d
  \ge 2k = n - 1\), as desired.
\end{proof}

Recall \(d(G)\) denotes the smallest \(d\) such that one can find a faithful
\(G \hookrightarrow \GL_d(\CC)\), and \(d(\Sigma) = d(\PMod(\Sigma))\). Let
\(\Sigma_{g, 1}\) be the closed genus \(g\) surface with a single marked point.
By replacing the number \(2\) by \(\min \{ d(E) : E \text{ is a cyclic
extension of } \Mod(\Sigma_{\lfloor \sfrac{g}{n} \rfloor, 1}) \}\) in the proof
of Proposition~\ref{thm:free-group-prod-linear-bound} we obtain
Theorem~\ref{thm:bigger-bounds}.

\begin{theorem}[Theorem~\ref{thm:bigger-bounds}]\label{thm:bigger-bounds-real}
  Let \(n \ge 1\) and \(g \ge 2n\). Then \(d(\Sigma_g^1) \ge n \cdot \min \{
  d(E) : E \text{ is a cyclic extension of } \Mod(\Sigma_{\lfloor
  \sfrac{g}{n} \rfloor, 1}) \}\).
\end{theorem}

\begin{proof}
  Let \(\rho : \Mod(\Sigma_g^1) \hookrightarrow \GL_d(\CC)\) be a faithful
  representation and \(g' = \lfloor \sfrac{g}{n} \rfloor\). Take
  \(d_{\operatorname{min}} = \min \{ d(E): E \text{ is a cyclic extension of }
  \Mod(\Sigma_{g', 1})\}\). We want to show \(d \ge n \cdot
  d_{\operatorname{min}}\). By passing to a smaller subsurface \(\Sigma_{n
  \cdot g'}^1 \subset \Sigma_g^1\) if necessary, we may assume \(g = n \cdot
  g'\) with \(g' \ge 2\). In that case, we may view \(\Sigma_g^1\) as an
  \((n+1)\)-hole sphere with \(n\) copies \(S_1, \ldots, S_n \subset
  \Sigma_g^1\) of \(\Sigma_{g'}^1\) attached along their boundaries.

  The natural maps \(\Mod(S_i) \to \Mod(\Sigma_g^1)\) are injective, so that we
  may regard \(\Mod(S_i)\) as a subgroup of \(\Mod(\Sigma_g^1)\). Since \(S_i
  \cap S_j = \emptyset\) for \(i \ne j\), there is a well-defined homomorphism
  \(\pi : \Mod(S_1) \times \cdots \times \Mod(S_n) \to \Mod(\Sigma_g^1)\) with
  \(\pi\!\restriction_{\Mod(S_j)} : \Mod(S_j) \to \Mod(\Sigma_g^1)\) injective
  for all \(i\).

  Take a maximal \((\Mod(S_1) \times \cdots \times \Mod(S_n))\)-invariant flag
  \begin{equation}\label{eq:bigger-bounds-invariant-flag}
    0 = W_0 \le W_1 \le \cdots \le W_{p+1} = \CC^d,
  \end{equation}
  so that the action \(\rho_i : \Mod(S_1) \times \cdots \times \Mod(S_n) \to
  \GL(W_{i+1}/W_i)\) of \(\Mod(S_1) \times \cdots \times \Mod(S_n)\) on each
  successive quotient is irreducible. Set \(W_{i+1}/W_i = W_{i,1} \otimes
  \cdots \otimes W_{i,n}\), where \(\rho_{i j} : \Mod(S_j) \to \GL(W_{ij})\) is
  an irreducible representation and \(\rho_i(f_1, \ldots, f_n) =
  \rho_{i,1}(f_1) \otimes \cdots \otimes \rho_{i,n}(f_n)\)
  \cite[Proposition~2.3.23]{kowalski}.

  For each \(j\), we may regard \(\rho_{ij}\) as a representation of
  \(\Mod(\Sigma_{g'}^1) \cong \Mod(S_j)\). Let us show that, for each \(j \le
  n\), we can find \(i\) with \(\dim W_{ij} \ge d_{\operatorname{min}}\). In
  that case,
  \[
    \begin{split}
      d
      & = \dim W_1/W_0 + \cdots + \dim W_{p+1}/W_p \\
      & \ge d_{\operatorname{min}}^{\#\{j:\dim W_{1,j} \ge d_{\operatorname{min}}\}}
          + \cdots
          + d_{\operatorname{min}}^{\#\{j:\dim W_{p,j} \ge d_{\operatorname{min}}\}} \\
      & \ge d_{\operatorname{min}}
          \cdot \# \{(i, j) : \dim W_{ij} \ge d_{\operatorname{min}}\} \\
      & \ge n \cdot d_{\operatorname{min}},
    \end{split}
  \]
  as desired.

  Fix \(j \le n\). Assume at first we can find \(i\) such that \(\ker \rho_{ij}
  \le \Mod(\Sigma_{g'}^1)\) is central. Since \(g' \ge 2\), the center of
  \(\Mod(\Sigma_{g'}^1)\) is generated by \(T_d\), the Dehn twist about the
  boundary curve \(d = \partial \Sigma_{g'}^1\). As in
  \textsection\ref{sec:dehn-twists}, the quotient \(\Mod(\Sigma_{g'}^1)/T_d
  \cong \Mod(\Sigma_{g',1})\) is the mapping class group of the closed
  genus \(g'\) surface with one marked point. In particular,
  \(\Mod(\Sigma_{g'}^1)\) is a central extension of \(\Mod(\Sigma_{g', 1})\) by
  \(\langle T_d \rangle \cong \ZZ\).

  If \(\rho_{ij}\) is faithful then \(\dim W_{ij} \ge d(\Sigma_{g'}^1) \ge
  d_{\operatorname{min}}\) by definition. We may thus assume \(\ker \rho_{ij}
  \ne 1\), in which case it is freely generated by a power \(T_d^{k_i}\) of
  \(T_d\) with \(k_i \ge 1\). In that case, \(\Mod(\Sigma_{g'}^1)/\ker
  \rho_{ij} = \Mod(\Sigma_{g'}^1)/T_d^{k_i}\) is a central extension of
  \(\Mod(\Sigma_{g', 1})\) by \(\ZZ/k_i\). In particular, \(\dim W_{ij} \ge
  d_{\operatorname{min}}\) once again.

  We are left to consider the case where, for some \(j\), \(\ker \rho_{ij}\) is
  \emph{not} central for all \(i\). Let us show this situation cannot happen.
  Denoting \(K = \ker \rho_{1,j} \cap \cdots \cap \ker \rho_{p,j}\), it is
  clear \(K\) acts on \(\CC^d\) by operators which are upper triangular with
  respect to a basis adapted to the flag from
  (\ref{eq:bigger-bounds-invariant-flag}). Hence its \(k^{\text{th}}\) derived
  subgroup \(K^{(k)}\) lies in \(\ker \rho\) for large enough \(k\). But \(K\)
  contains a free subgroup by
  Proposition~\ref{thm:normal-subgrps-contain-free}. This implies \(K^{(k)} \ne
  1\), which contradicts the assumption \(\rho\) is faithful. We are done.
\end{proof}

We now focus our attention on the faithful representations of
\(\PMod(\Sigma)\). A simple count shows that the maximal size of a family
\(a_1, \ldots, a_n, b_1, \ldots, b_n \subset \Sigma_g\) as in
(\ref{eq:intersection-rels}) whose pairwise geometric intersection numbers are
\(\le 2\) is \(2n = 3g-2\) or \(2n = 3g-3\), depending on whether \(g\) is even
or odd, respectively. As a consequence,
Proposition~\ref{thm:free-group-prod-linear-bound} thus recovers lower bounds
similar to Korkmaz'.

These families can be obtained by viewing \(\Sigma_g\) as a \(g\)-holed sphere
attached to \(1\)-holed tori \(H_1, \ldots, H_g \subset
\Sigma_g\). The \(g\) first pairs \((a_i, b_i)\) are taken as \(a_i, b_i
\subset H_i\) intersecting once. The remaining pairs can be obtained
by subdividing the \(g\)-holed sphere into \(g - 2\) pairs of pants and
combining them into \(4\)-holed spheres \(S_{1}, \ldots,
S_{\lfloor \sfrac{g-2}{2}\rfloor} \subset \Sigma_g\). We then choose
\(a_{g+i}, b_{g+i} \subset S_{i}\) intersecting twice.

To move beyond Korkmaz' bound of \(3g-2\) we instead use a different strategy.
We consider a family of curves \(a_1, \ldots, a_{3g-3}, b_1, \ldots, b_{3g-3}
\subset \Sigma\) where the curves \(a_i\) come from a certain pants
decomposition of \(\Sigma_g\), while the curves \(b_j\) are, in some sense,
``complementary'' to the curves \(a_i\). Unlike the curves in
(\ref{eq:intersection-rels}), the curve \(b_j\) is allowed to intersect \(b_i\)
twice for \(i \ne j\). See \textsection\ref{sec:main-result} for a proper
definition.

Take \(\rho : \PMod(\Sigma) \to \GL_d(\CC)\) with \(d\) small enough. In
\textsection\ref{sec:eigenvals} and \textsection\ref{sec:main-result} we will
show that, unless \(\rho\) kills the generators of \(\SIP_0(\Sigma)\),
the matrices \(M_i = \rho(T_{a_i}) - 1\) and \(N_j = \rho(T_{b_j}) - 1\)
satisfy the ``annihilation relations''
\begin{align}\label{eq:nm-rels}
  N_j M_i & = 0 \iff i \ne j    &
  M_j N_i & = 0 \iff i \ne j    &
  M_j M_i & = 0 \; \forall i, j.
\end{align}
We emphasize the relations \(N_j N_i = 0\) are not part of (\ref{eq:nm-rels}).

We now establish a simple lower bound for \(d\) such that we can find matrices
\(M_1, \ldots, M_n, N_1, \ldots, N_n \in M_d(\CC)\) satisfying
(\ref{eq:nm-rels}). This will be used in \textsection\ref{sec:main-result} to
show that \(\rho\) is, in fact, forced to kill \(\SIP_0(\Sigma)\).

\begin{lemma}\label{thm:naive-lemma}
  Let \(M_1, \ldots, M_n, N_1, \ldots, N_n \in M_d(\CC)\) be operators such
  that \(N_j M_i = 0\) if an only if \(i \ne j\). Then \(\dim \sum_i \range M_i
  \ge n\).
\end{lemma}

\begin{proof}
  We proceed by induction in \(n\). The base case \(n = 1\) is clear. Now
  suppose the theorem holds for a given \(n\) and let us show the same holds
  for \(n + 1\). Given \(M_1, \ldots, M_{n + 1}, N_1, \ldots, N_{n + 1}\) as
  above, it follows from the induction hypothesis that
  \[
    \dim \sum_{i \le n} \range M_i \ge n.
  \]
  We now claim one can find \(v \in \range M_{n + 1}\) with \(v \notin \sum_{i
  \le n} \range M_i\), so that \(\dim \sum_{i \le n+1} \range M_i \ge 1 + \dim
  \sum_{i \le n} \range M_i \ge n + 1\).

  Indeed, since \(N_{n+1} M_{n+1} \ne 0\), there is \(w \in \CC^n\) with
  \(N_{n+1} M_{n+1} w \ne 0\). On the other hand, \(N_{n+1} (M_1 w_1 + \cdots +
  M_n w_n) = N_{n+1} M_1 w_1 + \cdots N_{n+1} M_n w_n = 0\) for all \(w_1,
  \ldots, w_n \in \CC^n\). In other words, \(N_{n+1} \!\restriction_{\sum_{i
  \le n} \range M_i} = 0\) and thus \(v = M_{n+1} w \notin \sum_{i \le n}
  \range M_i\). This concludes the inductive step.
\end{proof}

\begin{proposition}\label{thm:dim-inf-bound}
  Let \(M_1, \ldots, M_n, N_1, \ldots, N_n \in M_d(\CC)\) be nonzero operators
  subject to the annihilation relations (\ref{eq:nm-rels}). Then \(\dim \sum_i
  (\range M_i + \range N_i) \ge 3n - d\). In particular, \(2d \ge 3n\).
\end{proposition}

\begin{proof}
  It is clear from Lemma~\ref{thm:naive-lemma} that \(\dim \sum_i \range M_i
  \ge n\) and \(\dim \sum_i \range N_i \ge n\). Let us show that \(\dim \left(
  \sum_i \range M_i \right) \cap \left( \sum_i \range N_i \right) \le d - n\),
  so that
  \[
    \begin{split}
      \dim \sum_i (\range M_i + \range N_i)
      & = \dim \sum_i \range M_i
        + \dim \sum_i \range N_i \\
      & - \dim
          \left( \sum_i \range M_i \right)
          \cap
          \left( \sum_i \range N_i \right) \\
      & \ge 2n - (d - n) \\
      & = 3n - d.
    \end{split}
  \]

  Since \(M_i M_j = 0\) for all \(i\) and \(j\), \(\left( \sum_i \range M_i
  \right) \cap \left( \sum_i \range N_i \right) \le \bigcap_i \ker
  M_i = \ker \Phi\), where \(\Phi = \bigoplus_i M_i : \CC^d \to
  \bigoplus_i \range M_i\). By the second relation in (\ref{eq:nm-rels}),
  we can find \(w_i \in \CC^n\) such that \(M_i N_i w_i \ne 0\). On the other
  hand, \(M_j N_i = 0\) for \(j \ne i\) and, in particular, \(M_j N_i w_i =
  0\). Hence \(\Phi(v_i) \ne 0\) lies in the copy of \(\range M_i\)
  inside of the codomain of \(\Phi\) for \(v_i = N_i w_i\).

  Choosing one such \(w_i\) for each \(i = 1, \ldots, n\) we get that the
  vectors \(\Phi(v_1), \ldots, \Phi(v_n)\) are linearly independent,
  so that \(\rank \Phi \ge n\). Hence \(\dim \ker \Phi \le d -
  n\), as desired.
\end{proof}

\section{Eigenspaces of $T_a$}
\label{sec:eigenvals}

Fix some \(\rho : \PMod(\Sigma) \to \GL_d(\CC)\) with \(d \le 4g - 4\). In this
section we study the \(1\)-eigenspace of \(L_a = \rho(T_a)\) for some
nonseparating \(a \subset \Sigma\). We establish a lower bound for the
dimension of the \(1\)-eigenspace of \(L_a\). 

As a consequence, we obtain the fact the matrices \(M_i = L_{a_i} - 1\) and
\(N_i = L_{b_i} - 1\) associated to the aforementioned family \(a_1, \ldots,
a_{3g-3}, \ldots, b_1, \ldots, b_{3g-3} \subset \Sigma\) satisfy the first two
relations in (\ref{eq:nm-rels}) -- see Corollary~\ref{thm:nil-product-lemma}.
This will be used in \textsection\ref{sec:main-result} to apply
Proposition~\ref{thm:dim-inf-bound} to the matrices \(M_i\) and \(N_j\) as
above.

Given \(a \subset \Sigma\), we denote the \(\lambda\)-eigenspace of \(L_a\) by
\(E_\lambda^a\). We also take \(E_{\lambda, k} = \ker (L_a - \lambda)^k\), so
that \(E_{\lambda, 1}^a = E_\lambda^a\) and \(E_{\lambda, d}^d\) is the
generalized \(\lambda\)-eigenspace of \(L_a\). Recall from
\textsection\ref{sec:dehn-twists} that the Dehn twists about nonseparating \(a, b
\in \Sigma\) are conjugate in \(\PMod(\Sigma)\). In particular, \(L_a \sim
L_b\). We may thus pass from one nonseparating curve to the next when
performing our analysis.

We will call \(\rho\) \emph{unipotent} if \(1\) is the only eigenvalue of
\(L_a\) for some (and hence all) nonseparating \(a \subset \Sigma\).
Establishing the \emph{unipotency} of low-dimensional representations is a
crucial step in the classification theorems of Korkmaz, Kasahara and
Kaufmann--Salter--Zhang--Zhong. This is summarized in the following
proposition.

\begin{proposition}[Kaufmann--Salter--Zhang--Zhong, Proposition~6.1 \cite{kaufmann-25}]\label{thm:only-eigenval-is-1}
  Let \(\Sigma\) be a surface of genus \(g \ge 4\) and \(\rho : \PMod(\Sigma)
  \to \GL_d(\CC)\) with \(d \le 4g - 3\). Given \(a \subset \Sigma\)
  nonseparating, the only eigenvalue of \(L_a\) is \(1\).
\end{proposition}

Building on the work of Kaufmann--Salter--Zhang--Zhong, we establish a lower
bound for the following dimension of the \(1\)-eigenspace of \(L_a\).

\begin{proposition}\label{thm:eigenval-multiplicity-quota}
  Let \(\Sigma\) be a genus \(g\) surface and \(\rho : \PMod(\Sigma) \to
  \GL_d(\CC)\) be nontrivial. Suppose either of the following conditions are
  met:
  \begin{enumerate}
    \item \(g \ge 4\) and \(d < 4g - 4\), or
    \item \(g \ge 7\) and \(d \le 4g - 4\).
  \end{enumerate}
  If \(a \subset \Sigma\) is nonseparating then \(\dim E^a_1 > 2g - 2\).
\end{proposition}

\begin{corollary}\label{thm:nil-product-lemma}
  Let \(\Sigma\) be a genus \(g\) surface and \(\rho : \PMod(\Sigma) \to
  \GL_d(\CC)\) satisfying either {\normalfont (1)} or {\normalfont (2)} from
  Proposition~\ref{thm:eigenval-multiplicity-quota}. Take \(a, b \subset
  \Sigma\) disjoint with \(a\) nonseparating. Then \((L_a - 1)(L_b - 1) = 0\).
\end{corollary}

\begin{proof}[Proof of Corollary~\ref{thm:nil-product-lemma}]
  Denote by \(\Sigma_a\) the surface obtained from \(\Sigma\) by cutting
  across \(a\).

  The result clearly holds for trivial \(\rho\). We may thus assume \(\rho\) is
  nontrivial. We can find a basis for \(\CC^d\) under which
  \[
    \rho(f) =
    \left(
    \begin{array}{c|c}
      \rho_1(f) & *             \\ \hline
              0 & \bar{\rho}(f)
    \end{array}
    \right)
  \]
  for all \(f \in \PMod(\Sigma_a)\), where the top-left and bottom-right
  blocks correspond to the action of \(f\) on \(E^a_1\) and
  \(\CC^d/E^a_1\), respectively.

  Now since \(\rho\) is nontrivial, \(\dim E_1^a > 2g - 2\) by
  Proposition~\ref{thm:eigenval-multiplicity-quota} and thus \(\dim \CC^d/E^a_1
  < 2g - 2\). It follows from Theorem~\ref{thm:korkmaz-main-thms} that
  \(\bar{\rho}(f) = 1\). Given any \(b \subset \Sigma\) disjoint from \(a\),
  \(b \subset \Sigma_a\) and so we may write
  \[
    L_b - 1 =
    \left(
    \begin{array}{c|c}
      L_b\!\restriction_{E^b_1} - 1 & * \\ \hline
                                  0 & 0
    \end{array}
    \right)
  \]
  in this basis. In particular, \(\range (L_b - 1) \le E^a_1 = \ker (L_a -
  1)\).
\end{proof}

We now review some results needed for the proof of
Proposition~\ref{thm:eigenval-multiplicity-quota}.

\begin{lemma}[Jordan inequalities]\label{thm:jordan-ineqs}
  Let \(A \in M_d(\CC)\) and \(\lambda \in \CC\). Consider the flag
  \[
    0
    =   E_{\lambda, 0}
    \le E_{\lambda, 1}
    \le E_{\lambda, 2}
    \le \cdots
    \le E_{\lambda, d},
  \]
  where \(E_{\lambda, k} = \ker (A - \lambda)^k\). Then \(\dim E_{\lambda,
  k+1}/E_{\lambda, k} \le \dim E_{\lambda, k}/E_{\lambda, k-1}\) for all \(k =
  0, 1, \ldots, d - 1\).
\end{lemma}

\begin{proof}
  The map \(A - \lambda\) induces an injective linear map \((A -
  \lambda) : E_{\lambda, k+1}/E_{\lambda, k} \hookrightarrow E_{\lambda,
  k}/E_{\lambda, k-1}\) for all \(k\).
\end{proof}

\begin{lemma}[Korkmaz, Lemma~4.3 \cite{korkmaz-23}]\label{thm:eigenspace-pair-is-invariant}
  Let \(\Sigma\) be a surface of genus \(g \ge 2\) and \(\rho : \PMod(\Sigma)
  \to \GL_d(\CC)\). Fix two nonseparating curves \(a, b \subset \Sigma\)
  intersecting at a single point and suppose \(E_\lambda^a = E_\lambda^b\).
  Then \(E_\lambda^a\) is a \(\PMod(\Sigma)\)-invariant subspace.
\end{lemma}

We are now ready to prove Proposition~\ref{thm:eigenval-multiplicity-quota}.

\begin{proof}[Proof of Proposition~\ref{thm:eigenval-multiplicity-quota}]
  Take some nontrivial \(\rho : \PMod(\Sigma) \to \GL_d(\CC)\) with \(d \le 4g
  - 4\) as above, where \(\Sigma\) is a surface of genus \(g \ge 4\). Let \(a
  \subset \Sigma\) be nonseparating. We know from
  Proposition~\ref{thm:only-eigenval-is-1} that \(1\) is the only eigenvalue of
  \(L_a\). Let \(\Sigma' \cong \Sigma_{g-1}^1\) be a subsurface as in
  Figure~\ref{fig:eigenvals-subsurf-def}.

  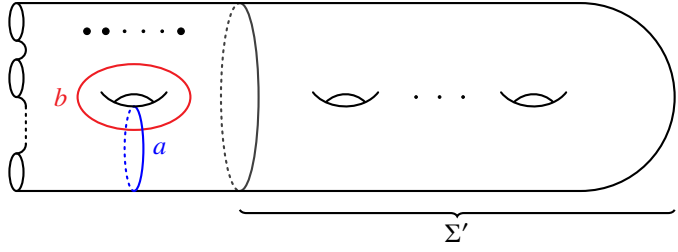
\begin{figure}
    \centering
    \begin{tikzpicture}[scale=1.25]
      \draw[darkgray, thick, line cap=round]
        ({2.25/2},1) arc[start angle=90, end angle=-90, y radius=1, x radius=.2];
      \draw[darkgray, thick, dotted, line cap=round]
        ({2.25/2},1) arc[start angle=90, end angle=270, y radius=1, x radius=.2];
      \draw[thick, decoration={brace, mirror, raise=0.25cm}, decorate]
        ({2.25/2},-1) -- (5.75,-1) node[pos=.5, anchor=north, yshift=-.3cm]{$\Sigma'$};
      \draw[thick] (-1.25,1) -- (4.75, 1) arc (90:-90:1) -- (-1.25,-1);
      \draw[thick, line cap=round] (-1.25,.6)  arc (90:-90:.1)
                                   (-1.25,0)   arc (90:0:.1)
                                   (-1.25,-.6) arc (-90:0:.1);
      \draw[dotted, thick, line cap=round] (-1.15,-.1) -- (-1.15,-.5);
      \draw[thick]
        (-1.25,1)
        arc[start angle=90, end angle=-270, x radius=.075, y radius=.2]
        (-1.25,.4)
        arc[start angle=90, end angle=-270, x radius=.075, y radius=.2]
        (-1.25,-.6)
        arc[start angle=90, end angle=-270, x radius=.075, y radius=.2];
      \draw[Red, thick]
        (-90:.35) arc[x radius=.6, y radius=.35, start angle=-90, end angle=270];
      \draw[Red] (180:.6) node[left] {$b$};
      \draw[blue, thick, line cap=round]         (0,-1)        arc[x radius=.1, y radius=.45, start angle=-90,  end angle=90];
      \draw[blue, thick, line cap=round, dotted] (0,-1)        arc[x radius=.1, y radius=.45, start angle=-90,  end angle=-270];
      \draw[blue]                          (.1,{-1+.45}) node[right] {$a$};
      \torushole{0}{0}
      \torushole{2.25}{0}
      \foreach \i in {0,1,2} { \filldraw (3,0)+({\i*.25}, 0) circle(0.4pt); }
      \torushole{4.25}{0}
      \filldraw (-.5,.7) circle(1pt)
                (-.3,.7) circle(1pt)
                (.5,.7)  circle(1pt);
      \foreach \i in {0,1,2} { \filldraw (-.1,.7)+({\i*.2}, 0) circle(0.4pt); }
    \end{tikzpicture}
    \caption{The subsurface $\Sigma' \cong \Sigma_{g-1}^1$.}
    \label{fig:eigenvals-subsurf-def}
  \end{figure}

  Suppose by contradiction \(\dim E_1^a \le 2g - 2\). First, assume \(\dim E_1^a
  < 2g - 2\) and consider the
  \(\Mod(\Sigma')\)-invariant flag
  \[
    0 \le E_1^a = E_{1, 1}^a
    \le E_{1, 2}^a
    \le \cdots
    \le E_{1, d}^a = \CC^d,
  \]
  where \(E_{1, k}^a = \ker (L_a - 1)^k\) as above. By the Jordan inequalities
  (Lemma~\ref{thm:jordan-ineqs}), \(\dim E_{1, k+1}^a / E_{1, k}^a \le \dim
  E_1^a < 2g - 2\) for all \(k\). It thus follows from the flag triviality
  criterion (Lemma~\ref{thm:korkmaz-triviality-criteria}) that the restriction
  of \(\rho\) to \(\Mod(\Sigma')\) is trivial. But then \(L_a \sim L_c = 1\)
  for any nonseparating \(c \subset \Sigma'\), contradicting the assumption
  \(\dim E_1^a < 2g - 2\).

  It remains to show \(\dim E_1^a \ne 2g - 2\). Assume by contradiction \(\dim
  E_1^a = 2g - 2\) and denote by \(\rho_1 : \Mod(\Sigma') \to \GL(E_1^a)\) and
  \(\bar\rho: \Mod(\Sigma') \to \GL(\CC^d/E_1^a)\) the actions of
  \(\Mod(\Sigma')\) on \(E_1^a\) and \(\CC^d/E_1^a\), respectively. It
  follows from Theorem~\ref{thm:korkmaz-main-thms} that \(\rho_1\) and \(\bar
  \rho\) are either trivial or conjugate to the symplectic representation
  \(\Psi : \Mod(\Sigma') \to \Sp_{2g - 2}(\ZZ) \le \GL(H_\CC)\), where
  \(H_\CC = H_1(\Sigma_{g-1}; \CC)\).

  We consider three separate cases.

  \noindent\textbf{Case 1.} Assume \(\rho_1\) is trivial.
  In this case, \(L_c\!\restriction_{E_1^a} = \rho_1(T_c) = 1\) and thus
  \(E_1^a \le E_1^c\) for all nonseparating \(c \subset \Sigma'\).
  Since \(L_a\) and \(L_c\) are conjugate, this implies \(E_1^a = E_1^c\).
  By the same token, \(E_1^b = E_1^c = E_1^a\) for \(b\) as in
  Figure~\ref{fig:eigenvals-subsurf-def}. Now
  Lemma~\ref{thm:eigenspace-pair-is-invariant} implies \(E_1^a = E_1^b\) is
  \(\PMod(\Sigma)\)-invariant.

  We abuse the notation and denote by \(\rho_1 : \PMod(\Sigma) \to
  \GL(E_1^a)\) and \(\bar \rho: \PMod(\Sigma) \to \GL(\CC^d/E_a^1)\)
  the actions of \(\PMod(\Sigma)\) on \(E_1^a\) and \(\CC^d/E_a^1\),
  respectively. In that case, \(\rho_1\) and \(\bar \rho\) are both trivial
  by Theorem~\ref{thm:korkmaz-main-thms}. The flag triviality criterion
  (Lemma~\ref{thm:korkmaz-triviality-criteria}) applied to the flag \(0
  \le E_1^a \le \CC^d\) thus implies \(\rho\) is trivial, a
  contradiction.

  \noindent\textbf{Case 2.} Assume \(\rho_1 \sim \Psi\) and \(\bar\rho\) is
  trivial. In this case, we can find a basis for \(\CC^d\) under which
  \[
    \rho(f) =
    \left(
    \begin{array}{c|c}
      \Psi(f) & \begin{matrix}c_1(f)&c_2(f)&\cdots&c_{d-2g+2}(f)\end{matrix}\\ \hline
            0 & 1
    \end{array}
    \right)
  \]
  for all \(f \in \Mod(\Sigma')\). It is not hard to check the maps \(c_k :
  \Mod(\Sigma') \to E_1^a \cong H_\CC\) are crossed
  homomorphisms.

  Now Lemma~\ref{thm:crossed-hom-characterization} implies that, given \(c
  \subset \Sigma'\) nonseparating, \(c_k(c) = \mu_k \cdot c\) for some
  \(\mu_k \in \CC\). By tweaking the above basis, we can find a
  second basis for \(\CC^d\) under which
  \[
    L_c =
    \left(
    \begin{array}{c|c}
      \Psi(T_c) & \begin{matrix} \mu \cdot c & 0 & \cdots & 0 \end{matrix} \\ \hline
              0 & 1
    \end{array}
    \right)
  \]
  for some \(\mu \in \CC\). Hence \(\operatorname{codim} E_1^a =
  \operatorname{codim} E_1^c \le 2\) by (\ref{eq:symplectic-rep-formula}),
  a contradiction for \(d > 2g\).

  We may thus assume \(d \le 2g\), in which case
  Theorem~\ref{thm:korkmaz-main-thms} says \(\rho\) is either trivial or
  conjugate to the symplectic representation \(\Psi : \Mod(\Sigma) \to
  \Sp_{2g}(\ZZ)\). The former contradicts the hypothesis \(\rho\) is
  nontrivial, so \(\rho \sim \Psi\). But then \(\dim E_1^a = 2g - 1\) by
  (\ref{eq:symplectic-rep-formula}), contradicting the assumption \(\dim
  E_1^a = 2g - 2\).

  \noindent\textbf{Case 3.} Finally, assume \(\rho_1 \sim \bar \rho \sim
  \Psi\). This last case is only possible if \(d = 4g - 4\), which is only
  relevant to our proof when \(g \ge 7\). We thus assume \(d = 4g - 4\) and
  \(g \ge 7\) from now on.

  We regard \(\CC^{4g-4}\) as a \(\ZZ[\Mod(\Sigma')]\)-module, where
  \(\ZZ[\Mod(\Sigma')]\) denotes the group ring of \(\Mod(\Sigma')\) and
  \(f \in \Mod(\Sigma')\) acts on \(\CC^d\) by \(\rho(f)\), as in
  \textsection\ref{sec:twisted-cohomology}. In this case, \(\CC^{4g-4}\) is an
  extension of \(H_\CC = H_1(\Sigma_{g-1}; \CC)\) by \(H_\CC\). This means
  \(\CC^{4g-4}\) fits into a short exact sequence of the form
  \begin{equation}\label{eq:short-exact-seq}
    \begin{tikzcd}
      0          \rar &
      H_\CC      \rar &
      \CC^{4g-4} \rar &
      H_\CC      \rar &
      0.
    \end{tikzcd}
  \end{equation}
  Such extensions are classified by the group
  \(\operatorname{Ext}^1_{\ZZ[\Mod(\Sigma')]}(H_\CC, H_\CC) =
  \operatorname{Ext}^1_{\ZZ[\Mod(\Sigma')]}(H_\ZZ, H_\ZZ) \otimes_\ZZ
  \CC\), where \(H_\ZZ = H_1(\Sigma_{g-1}; \ZZ)\).

  On the one hand, \(\operatorname{Ext}^1_{\ZZ[\Mod(\Sigma')]}(H_\ZZ, H_\ZZ) =
  H^1(\Mod(\Sigma'); \operatorname{Hom}_\ZZ(H_\ZZ, H_\ZZ))\)
  \cite[Proposition~2.2]{brown}. Here \(\ZZ[\Mod(\Sigma')]\) acts on
  \(\operatorname{Hom}_\ZZ(H_\ZZ, H_\ZZ)\) by \(f \cdot A = \Psi(f) \circ A
  \circ \Psi(f)^{-1}\) for all \(f \in \Mod(\Sigma')\) and \(A \in
  \operatorname{Hom}_\ZZ(H_\ZZ, H_\ZZ)\). Hence \(\operatorname{Hom}_\ZZ(H_\ZZ,
  H_\ZZ) \cong H_\ZZ^* \otimes_\ZZ H_\ZZ\) as \(\ZZ[\Mod(\Sigma')]\)-modules,
  where \(\ZZ[\Mod(\Sigma')]\) acts on \(H_\ZZ^* =
  \operatorname{Hom}_\ZZ(H_\ZZ, \ZZ)\) via \(f \cdot A = A \circ
  \Psi(f)^{-1}\).

  On the other hand, the intersection pairing \(\langle \, , \rangle : H_\ZZ
  \times H_\ZZ \to \ZZ\) induces a \(\ZZ[\Mod(\Sigma')]\)-module isomorphism
  \(H_\ZZ^* \cong H_\ZZ\), so that \(H_\ZZ^* \otimes_\ZZ H_\ZZ \cong
  H_\ZZ^{\otimes 2}\). Since \(g-1 \ge 6\),
  \(\operatorname{Ext}^1_{\ZZ[\Mod(\Sigma')]}(H_\ZZ, H_\ZZ) \cong
  H^1(\Mod(\Sigma'); H_\ZZ^{\otimes 2}) = 0\) by
  Corollary~\ref{thm:coh-tensor-square-vanishes}. This implies the sequence
  (\ref{eq:short-exact-seq}) splits.

  We can thus find a basis for \(\CC^{4g-4}\) under which
  \[
    \rho(f) =
    \left(
    \begin{array}{c|c}
      \Psi(f) & 0       \\ \hline
            0 & \Psi(f)
    \end{array}
    \right)
  \]
  for all \(f \in \Mod(\Sigma')\). Taking \(f = T_c\) for some nonseparating
  \(c \subset \Sigma'\), we can see \(\dim E_1^a = \dim E_1^c = 4g - 6\), a
  contradiction.
\end{proof}

\section{Lower Bounds for Faithful Representations of the Mapping Class Group}
\label{sec:main-result}

In this section we conclude our proof of
Theorem~\ref{thm:bound-on-faithful-mcg}. Let \(\Sigma\) be a genus \(g\)
surface, possible with boundary components and marked points. We embed
\(\Sigma\) in \(\Sigma_g\) by capping the boundary components with disks.

Given \(\rho : \PMod(\Sigma) \to \GL_d(\CC)\) with \(d \le 4g - 4\), our goal
is showing \(\SIP_0(\Sigma) \le \ker \rho\). As mentioned before, our strategy
is to apply Proposition~\ref{thm:dim-inf-bound} to the family of matrices \(M_i
= \rho(T_{a_i}) - 1\) and \(N_i = \rho(T_{b_i}) - 1\) associated with \(\rho\),
where \(a_1, \ldots, a_{3g-3}, b_1, \ldots, b_{3g-3} \subset \Sigma\) are
obtained from a certain pants decomposition of \(\Sigma_g\). We begin by
defining the curves \(a_i, b_j \subset \Sigma\).

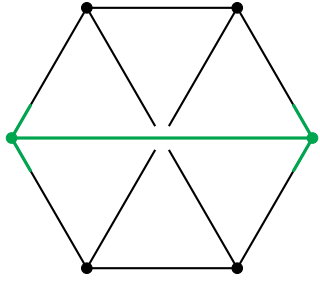
\begin{figure}
  \centering
  \begin{tikzpicture}
    \foreach \i in {1,2,3,4,5,6}
      \draw[thick] ({\i * (360/6)}:2) -- ({(\i + 1) * (360/6)}:2);
    \foreach \i in {1,2}
      \draw[thick] ({\i * (360/6)}:2) -- ({(\i + 3) * (360/6)}:2);
    \filldraw[white] circle (5pt);
    \draw ({3 * (360/6)}:2) -- ({6 * (360/6)}:2);
    \draw[Green, very thick, line cap=round]
      (0:2)   -- (180:2)
      (0:2)   -- +({-.5*cos(60)}, { .5*sin(60)})
      (0:2)   -- +({-.5*cos(60)}, {-.5*sin(60)})
      (180:2) -- +({ .5*cos(60)}, { .5*sin(60)})
      (180:2) -- +({ .5*cos(60)}, {-.5*sin(60)});
    \foreach \i in {1,2,4,5} \filldraw       ({\i * (360/6)}:2) circle (2pt);
    \foreach \i in {3,6}     \filldraw[Green] ({\i * (360/6)}:2) circle (2pt);
  \end{tikzpicture}
  \caption{The graph $\Gamma_g$ for $g = 4$, with the neighborhood $U_i$ of
  $e_i$ highlighted in green.}
  \label{fig:graph-def}
\end{figure}

Consider a trivalent graph \(\Gamma_g\) with \(2g - 2\) vertices given as
follows. We start by arranging \(2g - 2\) vertices uniformly in a circle, so
that, for each vertex we draw, we also draw its antipode. We then join each
vertex in the circle with the two adjacent vertices and its antipode, as in
Figure~\ref{fig:graph-def}. We embed \(\Gamma_g\) in \(3\)-space as to avoid
edge intersections.

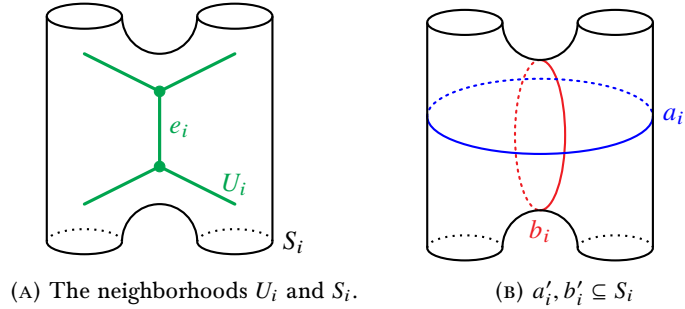
\begin{figure}
  \centering
  \begin{subfigure}[b]{.3\textwidth}
    \centering
    \begin{tikzpicture}
      \draw[thick] (0,  0) ellipse (.5 and {.5/3});
      \draw[thick] (2,  0) ellipse (.5 and {.5/3});
      \draw[thick] (.5, 0) arc     (180:360:.5);

      \draw[thick] (-.5, 0) -- +(0, -3)
                   (2.5, 0) -- +(0, -3) node[right]{$S_i$};

      \draw[Green, very thick, line cap=round]
        (0, -.5) -- (1, -1) -- (2, -.5)
        (1,  -2) -- (0, -2.5)
        (1,  -1) -- (1, -2)
        (1,  -2) -- (2, -2.5);

      \draw[Green] (1, -1.5) node[right]{$e_i$} (2, -2.5) node[above]{$U_i$};
      \filldraw[Green] (1, -1) circle(2pt) (1, -2) circle(2pt);

      \draw[thick]         (-.5, -3) arc[x radius=.5, y radius={.5/3}, start angle=180, end angle=360];
      \draw[thick, dotted] (-.5, -3) arc[x radius=.5, y radius={.5/3}, start angle=180, end angle=0];
      \draw[thick]         (1.5, -3) arc[x radius=.5, y radius={.5/3}, start angle=180, end angle=360];
      \draw[thick, dotted] (1.5, -3) arc[x radius=.5, y radius={.5/3}, start angle=180, end angle=0];
      \draw[thick]         (.5,  -3) arc (180:0:.5);
    \end{tikzpicture}
    \caption{The neighborhoods $U_i$ and $S_i$.}
  \end{subfigure}
  \begin{subfigure}[b]{.3\textwidth}
    \centering
    \begin{tikzpicture}
      \draw[Red, thick, line cap=round]         (1, -.5) arc[x radius={1/3}, y radius=1, start angle=90, end angle=-90]
                                                         node[below]{$b_i$};
      \draw[Red, thick, dotted, line cap=round] (1, -.5) arc[x radius={1/3}, y radius=1, start angle=90, end angle=270];
      \draw[blue, thick, line cap=round]         (-.5, -1.25) arc[x radius=1.5, y radius=.5, start angle=180, end angle=360] node[right]{$a_i$};
      \draw[blue, thick, dotted, line cap=round] (-.5, -1.25) arc[x radius=1.5, y radius=.5, start angle=180, end angle=0];
      \draw[thick] (0,  0) ellipse (.5 and {.5/3});
      \draw[thick] (2,  0) ellipse (.5 and {.5/3});
      \draw[thick] (.5, 0) arc     (180:360:.5);
      \draw[thick] (-.5, 0) -- +(0, -3) (2.5, 0) -- +(0, -3);
      \draw[thick]         (-.5, -3) arc[x radius=.5, y radius={.5/3}, start angle=180, end angle=360];
      \draw[thick, dotted] (-.5, -3) arc[x radius=.5, y radius={.5/3}, start angle=180, end angle=0];
      \draw[thick]         (1.5, -3) arc[x radius=.5, y radius={.5/3}, start angle=180, end angle=360];
      \draw[thick, dotted] (1.5, -3) arc[x radius=.5, y radius={.5/3}, start angle=180, end angle=0];
      \draw[thick]         (.5,  -3) arc (180:0:.5);
    \end{tikzpicture}
    \caption{$a_i, b_i \subset S_i$}
  \end{subfigure}
  \caption{Definition of the curves $a_i$ and $b_i$.}
  \label{fig:mcg-curves-def}
\end{figure}

The graph \(\Gamma_g\) is connected and has \(3g - 3\) edges \(e_1, e_2,
\ldots, e_{3g-3}\). What is more, for each such edge \(e_i\) we can find a
small ``double-Y-shaped'' neighborhood \(U_i \subset \Gamma_g\) of \(e_i\) such
that \(\Gamma_g \setminus U_i\) is still connected, as in
Figure~\ref{fig:mcg-curves-def}. By thickening \(\Gamma_g\) we obtain a
genus \(g\) handlebody with boundary \(\Sigma_g\).

For each \(e_i\), let \(a_i \subset \Sigma_g\) be a meridian around \(e_i\).
These curves form a pants decomposition of \(\Sigma_g\). By thickening the
neighborhood \(U_i\) of \(e_i\) we obtain a neighborhood \(S_i\) of \(a_i\)
which is a \(4\)-holed sphere. This neighborhood may also be obtained by gluing
two adjacent pairs of pants along their common boundary component \(a_i\). Let
\(b_i \subset \Sigma_g\) be as in Figure~\ref{fig:mcg-curves-def}.

We can find representatives of \(a_i\) and \(b_i\) lying on \(\Sigma \subset
\Sigma_g\) and avoiding its marked points, which we denote by \(a_i, b_i
\subset \Sigma\) as well. Although such curves depend on a choice of
representative of \(a_i, b_i \subset \Sigma_g\), this choice is inconsequential
to us. Since \(\Gamma_g \setminus U_i\) is connected, so is \(\Sigma_g
\setminus S_i\). In particular, \(a_i\) and \(b_i\) are nonseparating.

It is clear the curves \(a_i\) and \(a_j\) are disjoint for \(i \ne j\). While
\(a_i\) is disjoint from \(b_j\) for all \(i \ne j\), a careful reader might
notice the curves \(b_i\) and \(b_j\) \emph{can} intersect each other
nontrivially when \(e_i\) and \(e_j\) are adjacent edges of \(\Gamma_g\). We
will see this will not pose a problem to our proof.

We are now ready to prove Theorem~\ref{thm:bound-on-faithful-mcg}. Recall from
\textsection\ref{sec:torelli-johnson} that the subgroup \(\SIP_0(\Sigma)\)
generated by the simple intersection maps \([T_a, T_b]\) with \(\Sigma
\setminus (a \cup b)\) connected is a normal subgroup of
\(\mathcal{I}(\Sigma)\), normally generated by any such generator.

\begin{theorem}[Theorem~\ref{thm:bound-on-faithful-mcg}]\label{thm:bound-on-faithful-mcg-real}
  Let \(\rho : \PMod(\Sigma) \to \GL_d(\CC)\). If \(g \ge 4\) and \(d < 4g-4\)
  then \(\SIP_0(\Sigma_g) \le \ker \rho\). Moreover, if \(g \ge 7\) and \(d =
  4g - 4\) then \(\SIP_0(\Sigma_g) \le \ker \rho\).
\end{theorem}

\begin{proof}
  Given \(a \subset \Sigma_g\), denote \(L_a = \rho(T_a)\) and \(E_1^a = \ker
  (L_a - 1)\) as above. Suppose by contradiction \(\ker \rho\) does not contain
  \(\SIP_0(\Sigma)\). This means \(L_a\) and \(L_b\) do not commute for some
  (and hence all) nonseparating \(a, b \subset \Sigma\) intersecting twice with
  \(\Sigma \setminus (a \cup b)\) is connected.

  Take \(a_1, \ldots, a_{3g - 3}, b_1, \ldots, b_{3g - g} \subset \Sigma\) as
  above. We claim that the matrices
  \begin{align}\label{eq:ni-mi-def}
      M_i & = L_{a_i} - 1 & N_j & = L_{b_j} - 1
  \end{align}
  satisfy the conditions in (\ref{eq:nm-rels}). In that case, it follows from
  Proposition~\ref{thm:dim-inf-bound} that \(2d \ge 9g - 9\), a contradiction
  for \(d \le 4g - 4\). To establish the claim, notice \(\Gamma_g \setminus
  U_i\) is connected. This implies \(\Sigma_g \setminus S_i\) is connected.
  Hence so is \(\Sigma \setminus (a_i \cup b_i)\). In particular,
  \(\SIP_0(\Sigma)\) is normally generated by \([T_{a_i}, T_{b_i}]\) for any
  \(i\).

  It is clear from Corollary~\ref{thm:nil-product-lemma} that \(N_j M_i = M_j
  N_i = 0\) for \(i \ne j\) and \(M_j M_i = 0\) for all \(i, j\). Since the
  curves \(b_i\) and \(b_j\) can intersect each other nontrivially, we cannot
  conclude \(N_i N_j = 0\). This fact is, however, immaterial to our proof, as
  the relations \(N_i N_j = 0\) are \emph{not} part of the hypothesis of
  Proposition~\ref{thm:dim-inf-bound}.

  On the other hand, by assumption, \(L_{a_i}\) and \(L_{b_i}\) do not commute.
  Hence \(M_i\) and \(N_i\) don't commute. In particular, \(M_i N_i \ne 0\) or
  \(N_i M_i \ne 0\). But \(M_i N_i = 0 \iff N_i M_i = 0\). Indeed, the pairs
  \((T_{a_i}, T_{b_i})\) and \((T_{b_i}, T_{a_i})\) are conjugate in
  \(\PMod(\Sigma)\): we can find \(f \in \PMod(\Sigma)\) with \(f(a_i) = b_i\)
  and \(f(b_i) = a_i\), so that \(f T_{a_i} f^{-1} = T_{b_i}\) and \(f T_{b_i}
  f^{-1} = T_{a_i}\). The operators \(M_i N_i\) and \(N_i M_i\) are thus
  conjugated by \(\rho(f) \in \GL_d(\CC)\). We are done.
\end{proof}

\subsection*{Disclaimer}

Co-funded by the European Union. Views and opinions expressed are
however those of the author(s) only and do not necessarily reflect those of the
European Union. Neither the European Union nor the granting authority can be
held responsible for them.

\printbibliography

\end{document}